\documentclass{article}
\usepackage[%
journal=AIHPD,    
lang=british,   
]{ems-journal}
\usepackage{graphicx, hyperref}
\usepackage{listings}
\usepackage{xcolor}

\usepackage{amsmath}
\usepackage{newtxtext,newtxmath}
\usepackage{amsmath, amsthm}
\usepackage{verbatim}
\usepackage{svg}
\usepackage[utf8]{inputenc}
\usepackage{amsmath,amsthm, geometry, wasysym}
\usepackage{graphicx, verbatim, caption, subcaption}
\usepackage{xcolor, svg}
\usepackage{fullpage}
\usepackage{hyperref}
\usepackage{tikz}
\usetikzlibrary{arrows.meta, positioning}
\usepackage{graphicx}
\usepackage{subcaption}  
\usepackage[section]{placeins} 

\usepackage[psamsfonts]{eucal}
\usepackage{microtype}
\usepackage{mathtools}

\definecolor{codegreen}{rgb}{0,0.6,0}
\definecolor{codegray}{rgb}{0.5,0.5,0.5}
\definecolor{codepurple}{rgb}{0.58,0,0.82}
\definecolor{backcolour}{rgb}{0.95,0.95,0.92}

\lstdefinestyle{mystyle}{
    backgroundcolor=\color{backcolour},   
    commentstyle=\color{codegreen},
    keywordstyle=\color{magenta},
    numberstyle=\tiny\color{codegray},
    stringstyle=\color{codepurple},
    basicstyle=\ttfamily\footnotesize,
    breakatwhitespace=false,         
    breaklines=true,                 
    captionpos=b,                    
    keepspaces=true,                 
    numbers=left,                    
    numbersep=5pt,                  
    showspaces=false,                
    showstringspaces=false,
    showtabs=false,                  
    tabsize=2
}

\author{Aalok Gangopadhyay and Hariharan Narayanan\\School of Technology and Computer Science\\Tata Institute of Fundamental Research\\Mumbai, India}
\date{}

\newtheorem{theorem}{Theorem}
\newtheorem{lemma}[theorem]{Lemma}

\newcommand{\E}{\ensuremath{\mathbb E}}
\newcommand{\R}{\ensuremath{\mathbb R}}

\newcommand{\lab}{\label}  \newcommand{\ra}{\ensuremath{\rightarrow}}  \def\a{{\mathbf{\alpha}}} \def\de{{\mathbf{\delta}}} \def\De{{{\Delta}}}  
 \def\var{{\mathrm{var}}} \def\beq{\begin{eqnarray}} \def\eeq{\end{eqnarray}} \def\ben{\begin{enumerate}}
\def\een{\end{enumerate}} \def\bit{\begin{itemize}}
\def\bel{\begin{lemma}}
\def\eel{\end{lemma}}
\def\eit{\end{itemize}} \def\beqs{\begin{eqnarray*}} \def\eeqs{\end{eqnarray*}} \def\bel{\begin{lemma}} \def\eel{\end{lemma}}
\newcommand{\N}{\mathbb{N}} \newcommand{\Z}{\mathbb{Z}}   
\newcommand{\T}{\mathbb{T}}      \renewcommand{\b}{\mathbf{b}} 

 \newcommand{\I}{I}   \newcommand{\p}{\mathbb{P}}
\newcommand{\PP}{\mathcal P}    \newcommand{\one}{\mathrm{1}}
  \newcommand{\la}{\lambda}  
  \def\eps{{\epsilon}}  \def\ie{i.\,e.\,} \def\g{G}
\def\vol{\mathrm{vol}}

\renewcommand{\u}{u}

\renewcommand{\g}{\mathcal{G}}

\newcommand{\ent}{\mathrm{ent}}
\newcommand{\hess}{\nabla^2}

\newtheorem{proposition}{Proposition}
\newtheorem{thm}{Theorem}

\newtheorem{definition}{Definition}
\newtheorem{remark}{Remark}
\newtheorem{claim}{Claim}
\newtheorem{observation}{Observation}
\numberwithin{equation}{section}
\numberwithin{figure}{section}
\renewcommand{\a}{\alpha}
\renewcommand{\b}{\beta}
\renewcommand{\g}{\gamma}
\newcommand{\diag}{\mathrm{diag}}

\newcommand{\spec}{\mathrm{spec}}
\newcommand{\Spec}{\mathrm{Spec}}
\newcommand{\Leb}{\mathrm{Leb}}

\newcommand{\tla}{\tilde{\lambda}}
\newcommand{\tmu}{\tilde{\mu}}
\newcommand{\tnu}{\tilde{\nu}}
\newcommand{\tilh}{\tilde{h}}

\newcommand{\lacl}{ \lambda^{\mathrm{cl}}}

\newcommand{\rel}{\to}

\newcommand{\wt}{\mathbf{wt}}

\newcommand{\HIVE}{\mathtt{HIVE}}
\newcommand{\AHIVE}{\mathtt{AUGHIVE}}
\newcommand{\HORN}{\mathtt{HORN}}

\newcommand{\GT}{\mathtt{GT}}

\newcommand{\weight}{\mathbf{wt}}
\newcommand{\edge}{{\diamond}}
\renewcommand{\T}{T}

\numberwithin{equation}{section}

\begin{document}

\title{On the randomized Horn problem and the surface tension of hives}



\emsauthor{1}{
	\givenname{Aalok}
	\surname{Gangopadhyay}
	\orcid{0000-0002-1876-3140}}{Aalok}
\emsauthor{2}{
	\givenname{Hariharan}
	\surname{Narayanan}
	\mrid{796121}
	\orcid{0000-0003-0568-7350}}{Hari}

\Emsaffil{1}{
	\department{School of Technology and Computer Science}
	\organisation{TIFR}
	\rorid{https://ror.org/03ht1xw27}
	\address{https://ror.org/03ht1xw27}
	\zip{400005}
	\city{Mumbai}
	\country{India}
	\affemail{aalok@alumni.iitgn.ac.in}}
\Emsaffil{2}{
	\department{School of Technology and Computer Science}
	\organisation{TIFR}
	\rorid{https://ror.org/03ht1xw27}
	\address{https://ror.org/03ht1xw27}
	\zip{400005}
	\city{Mumbai}
	\country{India}
	\affemail{hariharan.narayanan@tifr.res.in}}


\begin{abstract}Given two nonincreasing $n$-tuples of real numbers $\la_n$, $\mu_n$, the Horn problem \cite{Horn} asks for a description 
    of all nonincreasing $n$-tuples of real numbers $\nu_n$ such that there exist Hermitian matrices $X_n$, $Y_n$ and 
    $Z_n$ respectively with these spectra such that $X_n + Y_n = Z_n$.  There is also a randomized version of this problem where $X_n$ and $Y_n$ are sampled uniformly at random from orbits of Hermitian matrices arising from the conjugacy action by elements of the unitary group. One then asks for a description of the probability measure of the spectrum of the sum $Z_n$. Both the original Horn problem and its randomized version have solutions using the hives introduced by Knutson and Tao \cite{CZ, KT1, KT2}. In an asymptotic sense, as $n \ra \infty$, large deviations for the randomized Horn problem were given in \cite{NarSheff} in terms of the surface tension of hives. In this paper, we provide upper and lower bounds on this surface tension function. We also obtain a closed-form expression for the minimum surface tension integral over all continuum hives with boundary conditions arising from the semicircle law. Finally, we give several empirical results for random hives and lozenge tilings arising from an application of the octahedron recurrence from \cite{NarSheffTao} for large $n$ and a numerical approximation of the surface tension function.
    \end{abstract}

\maketitle

\tableofcontents
\section{Introduction}
\label{sec:Introduction}

The study of Hermitian matrices and their spectra, in particular, the Horn problem, which asks for the possible spectra of a sum of two Hermitian matrices given their individual spectra, has connections to quantum mechanics, representation theory, and statistical physics. 

In the  randomized version of the Horn problem that we consider, the Hermitian matrices $X_n$ and $Y_n$ are sampled independently from  unitarily invariant distributions, and one asks for a description of the spectrum of $X_n + Y_n$. Large deviations for the spectra of large random matrices were shown in \cite{NarSheff} to be related to ideas of free energy and surface tension that arose originally in statistical physics, and certain random hives that arise naturally here were shown to concentrate in \cite{NarSheffTao} and the existence of a limit was recently established in \cite{Nar}.

In this paper, we investigate the surface tension function in the randomized Horn problem, providing both upper and lower bounds. We also derive a closed-form expression for the entropy of surface tension-minimizing continuum hives with boundary conditions that arise from scaling limits of GUE. We present numerical simulations of  random hives as well as some empirical results on their mean and variance. We also provide simulations of lozenge tilings that arise as weight maximizing dimer coverings of certain domains, where the weights are stochastic (arising from the interlacing gaps of GUE minor processes) which were used in \cite{NarSheffTao}. Finally, we provide a numerical approximation of the surface tension function. For further simulations see \cite{Gang}.

\subsection{Hives}
\lab{sec:prelim}




Suppose $\alpha, \beta$ are Lipschitz strongly concave functions\footnote{We assume strong concavity, \ie that there exists a positive $\de$ such that $\a + \de x^2 $ and $\beta + \de x^2$ are concave.}  from $[0, 1]$ to $\R$ and $\g$ is a concave function from $[0, 1]$ to $\R$, such that $\a(0) = \b(0) = \g(0) = \a(1) = \b(1) = \g(1) = 0.$ 
 \begin{definition}
 Let $U$ be a subset of $\R$.  We say that a function $f:U \ra \R$ is strongly decreasing if there exists $\de > 0$ such that $f(x) + \de x$ is monotonically decreasing.
 \end{definition}
For a $n \times n$ Hermitian matrix $W$, let $\spec(W)$ denote the vector in $\R^n$ whose coordinates are the eigenvalues of $W$ listed in non-increasing order. We denote by $\Spec_n$ the set of all vectors in $\R^n$ whose coordinates are non-increasing and sum to $0$.

\begin{definition}\lab{def:lan}
Let $\la = \partial^- \a$, $\mu = \partial^- \b$ on $(0, 1]$ and $\nu = \partial^- \g,$ at all points of $(0, 1]$, where $\partial^-$ is the left derivative. We note that $\la$ and $\mu$ are strongly decreasing and $\nu$ is monotonically decreasing.
Let $\la_n(i) := n^2\left(\a\left(\frac{i}{n}\right)-\a\left(\frac{i-1}{n}\right)\right)$, for $i \in [n]$, and similarly,
$\mu_n(i) := n^2\left(\b\left(\frac{i}{n}\right)-\b\left(\frac{i-1}{n}\right)\right)$, and 
$\nu_n(i) := n^2\left(\g\left(\frac{i}{n}\right)-\g\left(\frac{i-1}{n}\right)\right)$. 
\end{definition}

\begin{definition}
We define the seminorm $\|\cdot\|_\I$ on $\R^n$ as follows.
Given  $\la_n(1) \geq \dots \geq \la_n(n)$, such that $\sum_i \la_n(i) = 0$, let $\|\la_n\|_\I$ be defined to be the supremum over all $1 \leq i \leq n$ of $\la_n(1) + \dots + \la_n(i)$. Let $\one$ denote the vector in $\R^n$ whose coordinates are all $1$. Given any vector $v \in \R^n$ such that $v - (n^{-1}) \one \langle v, \one\rangle$ has coordinates that are a permutation of $\la_n(1), \dots, \la_n(n)$, we set $\|v\|_\I := \|\la_n\|_\I,$ and $\|\cdot\|_\I$ is readily seen to be a seminorm.
\end{definition}

\begin{definition}
We use $\|\cdot\|_\I$ to also denote a seminorm on the space of functions of bounded variation (\ie expressible as the difference of two monotonically decreasing functions) 
where the value on a given function $f:[0,1]\ra\R$  equals $\sup_{n \in \N} n^{-2}\|f_n\|_\I$, where $f_n$ is obtained from $f$ via the discretization process in Definition~\ref{def:lan}. Note that 
if $f$ is monotonically decreasing, then $\|f\|_\I = \lim_{n \ra \infty}n^{-2}\|f_n\|_\I$. Let $B_I(\nu, \eps):= \{\nu'\in \R^n | \|\nu' -\nu\|_I < \eps\}.$ 
\end{definition}

\label{subsec:Hives}

\begin{figure}
\begin{center}
\includegraphics[scale=0.5]{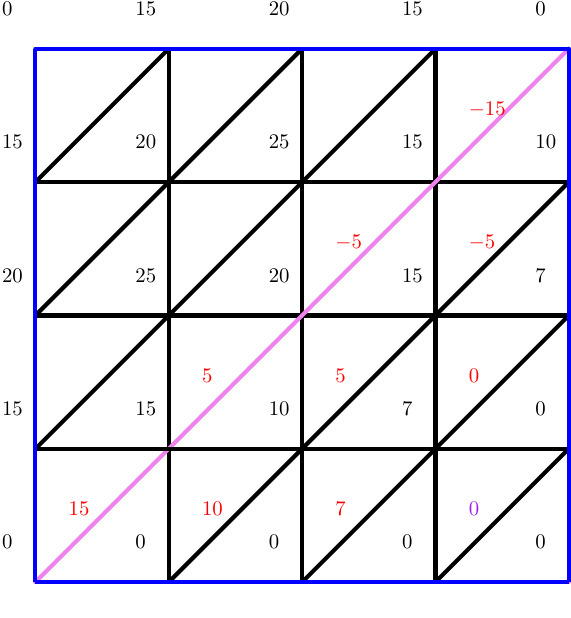}
\caption{An augmented hive for $n = 4$. The right triangle below the main diagonal of the $4\times 4$ square corresponds to a Gelfand-Tsetlin (GT) pattern (printed in red at the center of a diagonal edge) with top row $\nu_n$, given by the difference of the number
at a vertex and the one southwest to it.
The right triangle above the main diagonal corresponds to a hive with boundary conditions, $(\la_4, \mu_4, \nu_4)$.}
\label{fig:aug-hive}
\end{center}
\end{figure}

\begin{figure}[!h]
    \centering
    \includegraphics[width=0.5\linewidth]{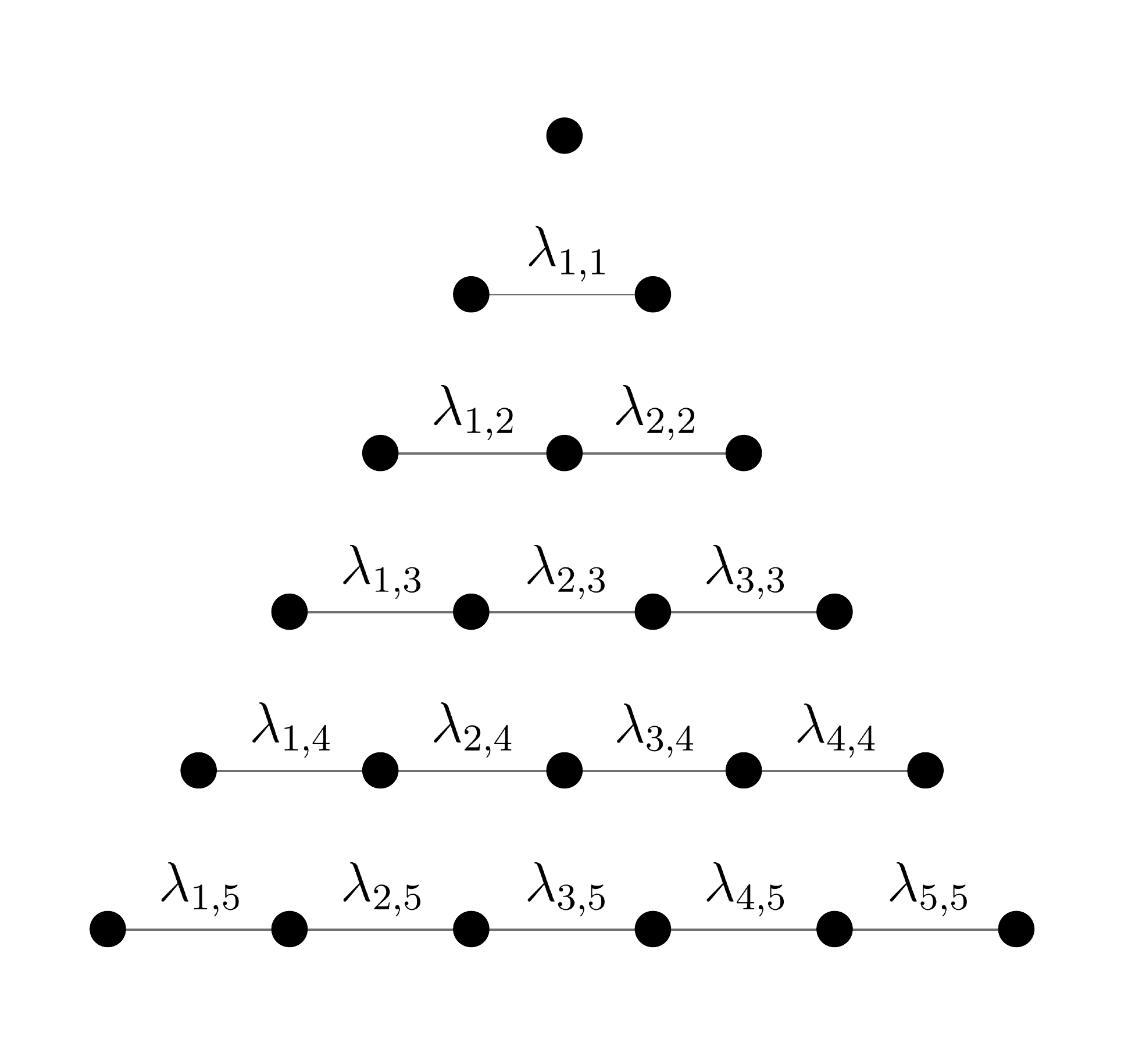}
    \caption{Gelfand-Tsetlin pattern for $n=5$.}
    \label{fig:GT_pattern}
\end{figure}

Let $T$ be the triangle $\{(x, y) \in [0, 1]^2| x \leq y\}.$ Let $T_n$ denote the set $nT \cap \Z^2.$
\begin{definition}[Discrete Hessian and the $\De_i$ on $\T_n$]
Let $f: \T_n \ra \R$ be a function.
\bit
\item Let $E_0(\T_n)$ be the set of all parallelograms $e_0\subseteq \T_n$ whose vertices are $\{(v_1, v_2),  (v_1 + 1, v_2),  (v_1 + 1, v_2 + 1),  (v_1 +2, v_2 + 1)\}.$

\item Let $E_1(\T_n)$ be the set of all parallelograms $e_1\subseteq \T_n$ whose vertices are $\{(v_1, v_2),  (v_1 + 1, v_2),  (v_1, v_2 + 1),  (v_1 +1, v_2 + 1)\}.$ 
\item Let $E_2(\T_n)$ be the set of all parallelograms $e_2\subseteq \T_n$ whose vertices are $\{(v_1, v_2),  (v_1 + 1, v_2+1),  (v_1, v_2 + 1),  (v_1 +1, v_2 + 2)\}.$ 
\eit

We define the discrete Hessian $\nabla^2(f):E(\T_n) \ra \R$ to be a  real-valued function on the set $E(\T_n) = E_0(\T_n) \cup E_1(\T_n) \cup E_2(\T_n)$ and the $\De_i$ from $\R^{\T_n}$ to $\R^{E_i(\T_n)}$ by 
\beq\lab{eq:A}
\nabla^2 f(e_0) := \De_0 f(e_0) :=  f(v_1, v_2) - f(v_1 + 1, v_2) - f(v_1 + 1, v_2 + 1) + f(v_1 +2, v_2 + 1).\nonumber\\
\nabla^2f(e_1) := \De_1 f(e_1) :=  - f(v_1, v_2) + f(v_1 + 1, v_2) + f(v_1, v_2 + 1) - f(v_1 + 1,  v_2 + 1).\nonumber\\
\nabla^2f(e_2) := \De_2 f(e_2) :=  f(v_1, v_2) - f(v_1+1, v_2+1) - f(v_1, v_2 + 1) + f(v_1 + 1, v_2 + 2).\nonumber\\
\eeq
\end{definition}

\begin{definition}[Rhombus concavity]\lab{def:rhomb}
Given a function $h:\T \ra \R$, and a positive integer $n$,  let $h_n$ denote the function from $\T_n$ to $\R$ such that for $(nx, ny) \in \T_n$, $h_n(nx, ny) = n^2 h(x, y).$
A function $h:\T \ra \R$ is called {\bf rhombus concave} if for any positive integer $n$, 
and any $i$, $\De_i h_n$ is nonpositive on $E_i(\T_n),$ and $h$ is continuous on $\T$.
The corresponding function $h_n$ is called {\bf discrete (rhombus) concave}.  
Note that a necessary and sufficient condition for a function $h_n$ from $\T_n$ to $\R$ to be discrete concave,  is that the piecewise linear extension (which we denote $\tilh_n$) of $h_n$ to $n\T$ is  concave.
Here each piece is an isosceles right triangle with a $\sqrt{2}-$length  hypotenuse parallel to the vector $(1, 1).$
\end{definition}

\begin{definition}[Hive]Let $ H_n(\lambda_n, \mu_n;  \nu_n)$ denote the set of all discrete concave functions $h_n:\T_n \ra \R,$ (which,  following Knutson and Tao \cite{KT1},  we call hives)
such that \ben \item $\forall i \in [n]\cup\{0\},\quad h_n(0,  i) = \sum_{j = 1}^i \la_n(j).$
\item $\forall i \in [n]\cup\{0\},\quad h_n(i,  n) = \sum_{j = 1}^n \la_n(j) + \sum_{j = 1}^i \mu_n(j).$
\item $\forall i \in [n]\cup\{0\},\quad h_n(i,  i) = \sum_{j = 1}^i \nu_n(j).$
\een
Let $|H_n(\lambda_n, \mu_n;  \nu_n)|$ denote the ${n-1 \choose 2}-$dimensional Lebesgue measure of this hive polytope. 
\end{definition}

\begin{definition}[$C^k$ hives] A $C^k$ hive for $k \geq 0$ is a rhombus concave map from $T$ to $\R$ that is the restriction of a $C^k$ function on some open neighborhood of $T$.
\end{definition}
\begin{definition}
A $C^0$ augmented hive is a function $a$ from $\Box:= [0,1]\times [0, 1]$, such that for each positive integer $n$, $a_n:\Box_n \ra \R$ defined by $a_n(nx,  ny) = n^2a(x,  y)$ for $(nx, ny) \in ([n]\cup\{0\}) \times ([n]\cup\{0\}),$  is a discrete augmented hive.
\end{definition}

\begin{definition}\lab{def:77old}
Let $\la = \partial^- \a$, $\mu = \partial^- \b$ on $(0, 1]$ be bounded, strongly decreasing functions and let $\nu = \partial^- \g,$ on $(0, 1]$ be a bounded, monotonically decreasing function.
We denote by $H(\la, \mu; \nu)$, the set of all $C^0$ hives from $T$ to $\R$, whose restrictions to the three sides are $\a, \beta$ and $\gamma$,  and set  $H(\la, \mu) := \cup_\nu H(\la, \mu; \nu)$. We denote the set of $C^0$ augmented hives with these boundary conditions by $A(\la, \mu; \nu).$
\end{definition}

\begin{definition}
    Let $V_n(\la_n(1), \dots, \la_n(n))$ denote the Vandermonde determinant  corresponding to the set of $n$ real numbers $\{\la_n(1), \dots, \la_n(n)\}$, whose value is $\prod\limits_{1 \leq i < j \leq n}(\la_n(i) - \la_n(j)).$ 
   \end{definition}
   \begin{definition}
   Let $\tau_n = \left(\frac{n-1}{2}, \frac{n-3}{2}, \dots, -\frac{n-3}{2}, -\frac{n-1}{2}\right).$
   \end{definition}
   
   Denoting probability densities with respect to the $n-1$ dimensional Lebesgue measure $$\Leb_{n-1, 0}(d\nu) = (d\nu(1)) \dots (d\nu(n-1))$$ on the hyperplane  in $\R^n$ consisting of vectors whose coordinates sum to $0$  by $\rho_n$, it is known through the work of Coquereaux and Zuber (see  Proposition 4 in \cite{CZ} and Equation (4) in \cite{Zuberhorn}, and also Knutson and Tao \cite{KT2} for a less explicit form of the result) that the following theorem holds.
   \begin{thm}[Coquereaux-Zuber]\lab{thm:1}
    Let $X_n = U_n \mathrm{diag}(\la_n)U_n^*$ and $Y_n = V_n \mathrm{diag}(\mu_n)V_n^*$ where $U_n$ and $V_n$ are independent random unitary matrices sampled from the Haar measure on the unitary group $\mathbb{U}_n.$ Then, 
   \begin{eqnarray}\lab{eq:2.4new} \rho_n\left[\spec(X_n + Y_n) = \nu_n\right] =  \frac{V_n(\nu_n)V_n(\tau_n)}{V_n(\la_n)V_n(\mu_n)} |H_n(\la_n, \mu_n; \nu_n)|.\end{eqnarray}
   \end{thm}

In Figure~\ref{fig:aug-hive}, a single point from the augmented hive polytope $A_4(\la_4, \mu_4; \nu_4)$ with boundary values $\la_4 = (15, 5, -5, -15)$, $\mu_4 = (15, 5, -5, -15)$, $\nu_4 = (15, 5, -5, -15)$ is depicted. Note that the rhombus inequalities corresponding to squares need not be satisfied at the pink diagonal or below it. Indeed they are violated at the square in the southeast corner.

   We view $$\bigcup\limits_{\nu'_n \in \Spec_n}A_n(\la_n, \mu_n; \nu'_n)$$
   as a $(n-1)^2 + (n-1)$ dimensional polytope contained in $\R^{(n+1)^2}$ with coordinates $(a_{ij})_{i, j \in [n]\cup \{0\}},$ corresponding to entries in an augmented hive, as in Figure~\ref{fig:aug-hive}. Thus, $j$ increases from bottom to top, and $i$ increases from left to right. Let $\Pi_{diag, n}:\R^{(n+1)^2} \ra \R^{n}$, denote the evaluation map that takes a $(n+1) \times (n+1)$  matrix $A$ to the first difference of its diagonal $(a_{i, i} - a_{i-1, i-1})_{1\leq i \leq n}.$ Let $\Pi_{diag}$ denote the continuous analogue of $\Pi_{diag, n},$ defined as follows. We view a continuum augmented hive as a function from $[0, 1]^2 \ra \R$.  Let $\Pi_{diag, n}:\R^{[0,1]^2} \ra \R^{(0, 1)}$, denote the evaluation map that takes a function $\tilde{A}$ from $[0, 1] \times [0, 1]$   to the left derivative on its diagonal (joining $(0, 1)$ to $(1, 1)$) $(\partial^-(\tilde A))$ in the direction $(1, 1)$.
   
   Recall that $\Spec_n$ denotes the cone of all possible vectors $x = (x_1, \dots, x_n)$ in $\R^n$ such that $$x_1 \geq x_2 \geq \dots \geq x_n,$$ and $\sum_i x_i = 0.$ 
   \begin{definition}Let $\eps_0 > 0$ and let $\la_n, \mu_n, \nu_n \in \Spec_n.$ Let $$G_n(\la_n, \mu_n; \nu_n, \eps_0 n^2) := \left(\bigcup\limits_{\nu'_n \in \Spec_n}A_n(\la_n, \mu_n; \nu'_n)\right)\cap \Pi_{diag, n}^{-1}(B_\I^n(\nu_n, \eps_0n^2)).$$
   Similarly, let $$G(\la, \mu; \nu, \eps_0) := \left(\bigcup\limits_{\nu'}A(\la, \mu; \nu')\right)\cap \Pi_{diag}^{-1}(B_\I(\nu, \eps_0)).$$
   \end{definition}
   
   In the remainder of this paper, by $o(1)$ we shall refer to a quantity whose magnitude tends to $0$ as $n \ra \infty.$

 \subsection{Surface tension of hives}
The theorem stated below in this section,  appeared in \cite{NarSheff}.  In what follows, we revisit the notation needed to state the theorem.
Let $\Leb_2$ denote twice the Lebesgue measure on $\R^2$ (there is an unfortunate typo in \cite{NarSheff} where $\Leb_2$ was defined to be the Lebesgue measure on $\R^2$).

Let $\la, \mu$ and $\nu$ be as in Definition~\ref{def:77old}.
There exists a convex function $\sigma$ that is termed the surface tension, defined on $\R_+^3$ (whose existence was proved in \cite{NarSheff}) and a functional $V$ defined in Lemma~\ref{lem:gamma} such that the following holds.
Here $h'$ is a continuum hive, and so $\hess h'$ is a Radon measure (see Definition 5 and Lemma 3 of \cite{NarSheff} for details). This Radon measure has a Lebesgue decomposition into a part that is absolutely continuous with respect to the Lebesgue measure $(\hess h')_{ac}$ and a singular part $(\hess h')_{sing}$ (see Section 4 of \cite{NarSheff}).
For a concave function $\g':[0, 1] \ra \R$ such that $\nu' = \partial^-\gamma'$,  we define the rate function 
\begin{eqnarray}\lab{eq:I} I(\g') :=  \ln\left(\frac{V(\la)V(\mu)}{V(\nu')}\right) + \inf\limits_{h' \in H(\la, \mu; \nu')}  \int_T \sigma((-1)(\hess h')_{ac})\Leb_2(dx).\end{eqnarray}

The following was proved in \cite{NarSheff}, with the notation of Definition~\ref{def:lan}.
\begin{thm}[Theorem 9, \cite{NarSheff}]\label{thm:4} Let $X_n, Y_n$ be independent random Hermitian matrices from unitarily invariant distributions with spectra $\la_n$,
$\mu_n$  respectively. 
Let $Z_n = X_n + Y_n$.
For any $\eps > 0$,
\begin{eqnarray}\lab{eq:master1}\lim\limits_{n \ra \infty}\left(\frac{-2}{n^2}\right)\ln \p_n\left[\spec(Z_n) \in {B}_\I^n(\nu_n, \eps)\right] = \inf\limits_{\partial^-\gamma'\in {B}_\I(\nu, \eps)}I(\g').\end{eqnarray}
\end{thm}

\subsection{Main results of this paper}
Set
   \[
   \Gamma(t) = \frac{2}{\pi} \int_{-1}^{t} \sqrt{1 - x^2} \, dx, \quad -1 \leq t \leq 1,
   \]
   and \( t = t(k, n) = \Gamma^{-1}(k/n) \) where \( k = k(n) \) is such that \( k/n \to a \in (0,1) \) as \( n \to \infty \). If \( x_k \) denotes the \( k \)-th eigenvalue in the GUE, then by a theorem of Gustavsson (see \cite[Theorem 1.1]{Gustavsson}), it holds that as \( n \to \infty \)
   \[
   \frac{x_k/\sqrt{2} - t \sqrt{2 n}}{\sqrt{\frac{\ln n}{4(1-t^2)n}}} \xrightarrow{d} N(0,1)
   \]
   in distribution. We have $x_k/\sqrt{2}$ in the above expression rather than $x_k$ as in \cite[Theorem 1.1]{Gustavsson}, because according to our convention, the standard deviation of the real and imaginary parts of a non-diagonal entry of a GUE matrix is $\frac{1}{\sqrt{2}}$ rather than $\frac{1}{2}$.
   
   \begin{definition}\lab{def:10}
    Let $\la^{\mathrm{cl}}:[0, 1] \ra \R$ be given by
   $\la^{\mathrm{cl}}(x) = 2\Gamma^{-1}(1 - x)$ for $x \in [0, 1]$. Note that the discrete counterpart as specified in Definition~\ref{def:lan} satisfies  $\la_n^{\mathrm{cl}}(i) = 2n^2 \int_{1 - \frac{i}{n}}^{1 - \frac{i-1}{n}}\Gamma^{-1}(x)dx$ for $i \in [n]$.
   \end{definition}

Let
   $\De (\sigma_\la, \sigma_\mu, \sigma_\nu)$ be the area of a two dimensional Euclidean triangle with the specified side lengths. In case there does not exist a triangle with the three arguments as the lengths of its sides, we set $\De(\sigma_\la, \sigma_\mu, \sigma_\nu) = 0$.
Let the functional $V$ be defined as in Lemma~\ref{lem:gamma}.  As shown in (\ref{eq:equality2}), $$\ln V(\sigma_\la \la^{\mathrm{cl}}) = \frac{5}{4} + \ln \sigma_\la.$$

The first main theorem of this paper is the following.

\begin{thm}\lab{thm:6}
 Let $\sigma_\la, \sigma_\mu, \sigma_\nu \in \R$, such that there exists a Euclidean triangle with these sidelengths and  $ \De  (\sigma_\la, \sigma_\mu, \sigma_\nu) > 0$.  Let $\la, \mu$ and $\nu$ denote bounded,  strongly decreasing  functions from $[0, 1]$ to $\R$.  Suppose that $\|\la\|_{L^2} = \sigma_\la$, $\|\mu\|_{L^2} = \sigma_\mu$ and $\|\nu\|_{L^2} = \sigma_\nu$. 

Then,
  \beq  \lab{eq:main1}  \inf\limits_{h \in H(\la, \mu; \nu)}  \int_T \sigma((-1)(\hess h)_{ac})\Leb_2(dx) & \geq & \nonumber\\ \ln V(\la) + \ln V(\mu) + \ln V(\nu)  - 5 &-& \ln\left(4  \De  (\sigma_\la, \sigma_\mu, \sigma_\nu)^2 \right).\eeq
 If
   $\la, \mu, \nu$ are respectively $\sigma_\la \la^{\mathrm{cl}}, \sigma_\mu \la^{\mathrm{cl}}, \sigma_\nu \la^\mathrm{cl}$, then equality holds.
   \end{thm}
Let $\varUpsilon$ be the nonnegative function on $\R^3_+$ defined in Section~\ref{sec:bounds_surf}.
The second  main result of the paper is:

\begin{thm}\lab{thm:9} For $s=(s_0,s_1,s_2) \in \R_+^3$, we have
\begin{eqnarray*}  \exp(\frac{5}{4}-\varUpsilon)\frac{(s_0 + s_1 + s_2)s_0s_1s_2}{(s_0 + s_1)(s_1 + s_2)(s_2 + s_0)} \leq  \exp(-\sigma(s))  \leq    \frac{\exp(5)(s_0 + s_1 + s_2)s_0s_1s_2}{36(s_0 + s_1)(s_1+s_2)(s_2 + s_0)}.\end{eqnarray*}
\end{thm}

\subsection{Broad dependency map of lemmas and theorems}

The diagram below summarizes how the core ingredients of the paper depend on one another.
Background large–deviation inputs inform the overall framework, while Sections~\ref{sec:Gaussian}--\ref{sec:bounds_surf} provide
the main technical steps leading to the two principal results.

\medskip

\begin{center}
\resizebox{\textwidth}{!}{%
\begin{tikzpicture}[
  >=Latex,
  node distance=12mm and 16mm,
  every node/.style={font=\small},
  res/.style={draw, very thick, rectangle, rounded corners=2pt, align=center, inner sep=3pt, fill=white},
  aux/.style={draw, thick, rounded corners=10pt, align=center, inner sep=3pt, fill=white},
  bgd/.style={draw, dashed, align=center, inner sep=3pt, fill=gray!10},
  arrow/.style={->, very thick}
]

\node[bgd] (T1) {Theorem~\ref{thm:4} (Narayanan--Sheffield)\\Large deviations / rate function};
\node[aux, right=of T1] (Vfun) {Definition of \(V(\cdot)\),\\ rate functional \(I(\cdot)\)};

\draw[arrow] (T1) -- (Vfun);

\node[aux, below left=20mm and -5mm of Vfun] (L1) {Lemma~\ref{lem:1}\\(expectations with negative variance)};
\node[res, right=20mm of L1] (Th5) {Theorem~\ref{thm:gaussian}\\(Gaussian integral identity)};
\node[aux, right=20mm of Th5] (L2) {Lemma~\ref{lem:Boltz}\\(Boltzmann / MaxEnt lemma)};
\node[res, right=20mm of L2] (Th6) {Theorem~\ref{thm:gaussian1}\\(Max entropy characterisation)};

\draw[arrow] (L1) -- (Th5);
\draw[arrow] (L2) -- (Th6);

\node[aux, below=22mm of Th5] (L3) {Lemma~\ref{lem:gamma}\\(limit \(\ln V_n(\lambda)-\ln V_n(\tau_n)\))};
\node[res, right=25mm of L3] (Th7) {Theorem~\ref{thm:gaussian2} \& \ref{thm:gaussian3} \& \ref{thm:5}\\(log-volume bounds for\\ hive polytopes)};

\draw[arrow] (Th5) -- (Th7);
\draw[arrow] (Th6) -- (Th7);
\draw[arrow] (L3) -- (Th7);

\node[res, below left=18mm and -6mm of Th7] (Main2) {Theorem~\ref{thm:6}\\(lower bound and equality case\\ via triangle area \(\Delta\))};
\node[res, right=28mm of Main2] (Main3) {Theorem~\ref{thm:9}\\(two-sided bounds on \(\sigma(s)\))};

\draw[arrow] (Th7) -- (Main2);
\draw[arrow] (L3) -- (Main2);
\draw[arrow] (Th7) -- (Main3);

\draw[arrow, dashed] (Vfun) |- (Th5);
\draw[arrow, dashed] (Vfun) |- (Th6);

\node[above=2mm of L1] {\textbf{Section~\ref{sec:Gaussian}}};
\node[above=2mm of L3] {\textbf{Section~\ref{sec:lim-log-vol}}};
\node[above=2mm of Main2] {\textbf{Section~\ref{sec:bounds_surf}}};
\end{tikzpicture}
}
\end{center}

\paragraph{Organisation of this paper.}
\begin{enumerate}
  \item Section~\ref{sec:Gaussian} establishes a Gaussian–integral identity (Theorem~\ref{thm:gaussian}), with Lemma~\ref{lem:1} handling the
        case of a negative variance parameter. In parallel, Lemma~\ref{lem:Boltz} yields an exponential-family
        maximum–entropy principle, culminating in Theorem~\ref{thm:gaussian1}.
  \item Section~\ref{sec:lim-log-vol} identifies the limiting Vandermonde functional \(V(\cdot)\) (Lemma~\ref{lem:gamma}) and proves a log–volume upper bound for hive polytopes (Theorem~\ref{thm:gaussian2}), a rigidity-box lower bound for the semicircle boundary (Theorem~\ref{thm:gaussian3}), and the resulting asymptotic equality in Theorem~\ref{thm:5}.  The older expected-GUE/Fradelizi argument is recorded separately for the positive raw-parameter regime.
  \item Section~\ref{sec:bounds_surf} uses Theorem~\ref{thm:gaussian2} (and Lemma~\ref{lem:gamma} where needed) to derive the main surface–tension results:
        the sharp lower bound with equality in the semicircle case (Theorem~\ref{thm:6}) and the general two–sided
        bounds on \(\sigma\) (Theorem~\ref{thm:9}).

\item In Section~\ref{recurrence-sec}, we present numerical experiments on random hives and the corresponding
lozenge tilings obtained via the octahedron recurrence.  These simulations use pairs of
independent GUE minor processes to generate random boundary data for hives, from which
large-scale lozenge tilings are produced.  The pointwise mean and variance of these random
hives are computed for two ensembles --- GUE matrices of different variances
and random projections. 
\item  Section~\ref{sec:EstimatingSurfaceTension} provides a detailed numerical study of the surface–tension function $\sigma$ appearing in Theorem~\ref{thm:6}. The approach involves making a guess of the shape of surface tension minimizing continuum hive, that is (non-rigorously) justified by a sampling experiment in Figure~\ref{fig:GUE_3,4,5} and using this shape to numerically solve the identities appearing in the equality cases of Theorem~\ref{thm:6}.  \end{enumerate}

  \section{A Gaussian integral}\lab{sec:Gaussian}

   For future use, we note the following theorem of Fradelizi \cite{Fradelizi}. 
   \begin{thm}[Fradelizi]\lab{thm:frad}
    The density at the center of mass of a logconcave density on $\R^{n}$ is greater or equal to $e^{- n}$ multiplied by the supremum of the density over $\R^n$.
   \end{thm}

  We identify two cases, which however are handled in a unified fashion: \ben \item  $a, b, c > 0$, \item $a, b > 0$ and $-c^2 > a^2 + b^2$.\een
  Thus, in the second case, $c$ is imaginary.
\begin{thm}\lab{thm:gaussian} 
For $\la, \mu, \nu \in \Spec_n + \R\one$, (where $\one$ denotes the vector of all ones), let $$F(\la, \mu, \nu) := \frac{V_n(\lambda)V_n(\mu)V_n(\nu)}{V_n(\tau_n)} |H_n(\lambda, \mu; \nu)| \exp\left(-\frac{1}{2}\left(\frac{|\lambda|^2}{a^2} + \frac{|\mu|^2}{b^2} + \frac{|\nu|^2}{c^2}\right)\right),$$ and otherwise let $F(\la, \mu, \nu) = 0.$ 
Then,
\begin{equation*}
\int\limits_{\substack{\R^n \times \R^n \times\R^n\\\sum \la + \sum \mu = \sum \nu}} F(\la, \mu, \nu) \Leb(d(\la, \mu, \nu)) = (2\pi)^{n}\left(\frac{a^2b^2c^2}{a^2 + b^2 + c^2}\right)^{n^2/2}.   \end{equation*}
\end{thm}
\begin{proof}
For $a > 0$, let us first consider the GUE joint density for the eigenvalues, \cite[Equation 2.4.23]{AndersonGuionnetZeitouni} namely  $$GUE(\lambda, a) := \left(\frac{1}{(2\pi)^\frac{n}{2}}\right) a^{-n^2} \frac{V_n(\la)^2}{V_n(\tau_n)}\exp\left(-\frac{|\la|^2}{2 a^2}\right)$$ for $\la \in \Spec_n + \R\one$ and $0$ otherwise.


By Theorem~\ref{thm:1} we see that 
$$\rho_n(\nu|\la, \mu) = \frac{V_n(\nu)V_n(\tau_n)}{V_n(\la)V_n(\mu)} |H_n(\la, \mu; \nu)|.$$
Denoting $\left(\frac{1}{(2\pi)^\frac{n}{2}}\right) |a|^{-n^2}\left(\frac{1}{(2\pi)^\frac{n}{2}}\right) |b|^{-n^2}$ by $\tilde{C}$, we observe that

$$\int_{\R^n}\int_{\R^n} \tilde{C} F(\la, \mu, \nu) d\la d\mu   =$$ 

\beq \lab{eq:tildeC} \int_{\R^n}\int_{\R^n} GUE(\lambda, a)GUE(\mu, b)\rho_n(\nu|\la, \mu)\exp\left(-\frac{|\nu|^2}{2c^2}\right) d\la d\mu .\eeq
Note that for any fixed $\la, \mu$, the function $ \rho_n(\nu|\la, \mu)$ is a singular probability measure supported on the hyperplane $\sum \la + \sum \mu = \sum \nu.$
Consider a multivariate Gaussian random variable $(X, Y)$ supported on $\R^{m + m}$ where  $X = (X_1, \dots, X_m)$ and $Y = (Y_1, \dots, Y_m)$, and the collection $\{X_1, \dots, X_m, Y_1, \dots Y_m\}$ is jointly independent (and of course, jointly Gaussian). Then, as is well known, $X + Y$ is also a multivariate Gaussian and  $$\var(X_i + Y_i) = \var X_i + \var Y_i.$$
It follows (by considering the diagonal and upper triangular parts) that  the sum of two independent scaled GUE matrices is a scaled GUE matrix, and that the variance of this sum is the sum of the variances of the scaled GUE matrices that are the summands. Therefore, the above expression equals
\begin{eqnarray*} GUE\left(\nu, \sqrt{a^2 + b^2}\right)\exp\left( - \frac{|\nu|^2}{2c^2}\right) = \left(\frac{c}{\sqrt{a^2 + b^2 + c^2}}\right)^{n^2} GUE\left(\nu, \frac{c\sqrt{a^2 + b^2}}{\sqrt{a^2 + b^2 + c^2}}\right).\end{eqnarray*} 
 Finally, since $GUE\left(\nu, \frac{c\sqrt{a^2 + b^2}}{\sqrt{a^2 + b^2 + c^2}}\right)$ is a probability density, we see that 
 \begin{eqnarray} \int_{\R^n}\left(\frac{c}{\sqrt{a^2 + b^2 + c^2}}\right)^{n^2} GUE\left(\nu, \frac{c\sqrt{a^2 + b^2}}{\sqrt{a^2 + b^2 + c^2}}\right)d\nu & = & \nonumber\\\left(\frac{c}{\sqrt{a^2 + b^2 + c^2}}\right)^{n^2}.\lab{eq:last3}\end{eqnarray} In the above array of equations, the branch of the square root is chosen such that $\frac{c}{\sqrt{a^2 + b^2 + c^2 }}   > 0.$
The theorem follows.
\end{proof}

\subsection{Maximum entropy formulation}

Consider the probability density $\rho_{GI}(\la, \mu, \nu)$ on $$\{(\la, \mu, \nu) \in \R^n \times \R^n \times \R^n | \sum \la + \sum \mu = \sum \nu\}$$ whose density is given by
\beq \lab{eq:GI} \rho_{GI}(\la, \mu, \nu) = \left((2\pi)^{n} \left(\frac{a^2b^2c^2}{a^2 + b^2 + c^2}\right)^{n^2/2}\right)^{-1} F(\la, \mu, \nu).\eeq
By (\ref{eq:last3}) we see that the marginal $\int \rho_{GI}(\la, \mu, \nu)d\la d\mu$ on $\nu$ is proportional to $GUE\left(\nu, \sqrt{\frac{c^2(a^2 + b^2)}{a^2 + b^2 + c^2}}\right)$. It follows that $$\E |\nu|^2 = \frac{c^2(a^2 + b^2)n^2}{(a^2 + b^2 + c^2)},$$ where the expectation is with respect to the probability density $\rho_{GI}(\la, \mu, \nu)$.
Similarly, with respect to the same measure, if $c$ is real, we have $$\E |\la|^2 = \frac{a^2(b^2 + c^2)n^2}{(a^2 + b^2 + c^2)}, \quad \E |\mu|^2 = \frac{b^2(c^2 + a^2)n^2}{(a^2 + b^2 + c^2)},$$ but when $c$ is imaginary, this requires a separate proof, which appears below.
\begin{lemma}\lab{lem:1}
    Suppose $a, b > 0$ and $-c^2 > a^2 + b^2.$ Then, with respect to the probability density $\rho_{GI}(\la, \mu, \nu)$ defined in (\ref{eq:GI}), we have $$\E |\la|^2 = \frac{a^2(b^2 + c^2)n^2}{(a^2 + b^2 + c^2)}, \quad \E |\mu|^2 = \frac{b^2(c^2 + a^2)n^2}{(a^2 + b^2 + c^2)}.$$
\end{lemma}
\begin{proof}
Let $f, g$ be functions on the space of $n \times n$ Hermitian matrices given by $$f(X) = 2^{-n/2}(\pi)^{-\frac{n^2}{2}} a^{-n^2}\exp\left(-\frac{1}{2a^2}\|X\|_{HS}^2\right)$$ and 
$$g(X) = 2^{-n/2}(\pi)^{-\frac{n^2}{2}} |c|^{-n^2}\exp\left(-\frac{1}{2c^2}\|X\|_{HS}^2\right)$$ respectively. 
The spectral measure associated with $f$ is the GUE spectral measure with variance $a^2$ and that associated with $g$ is an infinite measure with spectral density $$GUE_-(\lambda, c) := \left(\frac{1}{(2\pi)^\frac{n}{2}}\right) |c|^{-n^2} \frac{V_n(\la)^2}{V_n(\tau_n)}\exp\left(\frac{|\la|^2}{2 |c|^2}\right)$$ for $\la \in \Spec_n + \R\one$ and $0$ otherwise.
This can be seen from the fact that $GUE_-$ is the pushforward of a measure that is constant on the orbits of Hermitian matrices under the the unitary group.
Then the convolution of $f$ and $g$ is a function $f \ast g$ on the space of $n \times n$ Hermitian matrices given by $f \ast g(X) = \int f(Y)g(X - Y)dY.$
By completing the square, we observe that when $-c^2 > a^2 + b^2$, the function $f\ast g$ is given by $$f \ast g(X) = 2^{-n/2}(\pi)^{-\frac{n^2}{2}} |c^2 + a^2|^{-\frac{n^2}{2}}\exp\left(- \frac{1}{ 2(a^2 + c^2)}\|X\|_{HS}^2\right).$$ 
Therefore, the spectrum associated with $f\ast g$  is associated with  the $GUE_-$  measure with parameter $\sqrt{ a^2 + c^2},$ which is an infinite measure. In order to get the law of $\mu$, as in (\ref{eq:tildeC}), we need to weight this measure by $\exp\left(-|\mu|^2/2b^2\right)$. 

Indeed, recalling the identity for the conditional density of $\mu$ given $(\nu, -\lambda)$,
\[
\rho_n(\mu \mid \nu, -\lambda)
= \frac{V_n(\mu)\,V_n(\tau_n)}{V_n(-\lambda)\,V_n(\nu)}\,\bigl|H_n(\nu,-\lambda;\mu)\bigr| = \frac{V_n(\mu)\,V_n(\tau_n)}{V_n(\lambda)\,V_n(\nu)}\,\bigl|H_n(\lambda,\mu; \nu)\bigr|,
\]
 denoting $\left(\frac{1}{(2\pi)^\frac{n}{2}}\right) |a|^{-n^2}\left(\frac{1}{(2\pi)^\frac{n}{2}}\right) |c|^{-n^2}$ by $\tilde{B}$, we observe that
\[
\int_{\R^n}\int_{\R^n} \tilde{B}\, F(\lambda,\mu,\nu)\, d\lambda\, d\nu
= \int_{\R^n}\int_{\R^n}
GUE(\lambda,a)\,GUE_-(\nu,c)\,\rho_n(\mu\mid\nu, - \lambda)\,
\exp\!\left(-\frac{|\mu|^2}{2b^2}\right)\, d\lambda\, d\nu .
\]

In view of the fact that $a^2 + b^2 + c^2 < 0$, we obtain the result that $\mu$ obeys the law of the GUE spectral measure with variance $\frac{b^2(a^2 + c^2)n^2}{a^2 + b^2 + c^2}$ as desired. The proof for $\la$ is similar.
\end{proof}
We now consider the cone of {\it triply augmented hives} $\mathbb{A}^3$ defined below. 
\begin{definition}\lab{def:triple}
 Let $$\mathbb{A}^3(\la, \mu; \nu) := GT(\la) \times GT(\mu) \times GT(\nu) \times H_n(\la, \mu; \nu).$$ where $\la, \mu, \nu \in \Spec_n + \R \one.$ By $\mathbb{A}^3,$  we denote the cone $\bigcup_{\la, \mu, \nu} \mathbb{A}^3(\la, \mu; \nu)$ where the union is taken over  $\la, \mu, \nu \in \Spec_n + \R \one$. We denote by $\phi:\mathbb{A}^3 \ra\{(\la, \mu, \nu): \sum \la + \sum \mu = \sum \nu\}$ the projection map that takes a triply augmented hive to its boundary values.
\end{definition}
 
 It is known (see for example \cite{NarSheffTao}) that the volume of $GT(\la)$ equals $V_n(\la)/V_n(\tau_n).$
\begin{definition}\label{def:p}
Let the probability distribution $p$ on $\mathbb{A}^3$ (with respect to the Lebesgue measure) at any point $(g_\la, g_\mu, g_\nu, h) \in \mathbb{A}^3(\la, \mu; \nu),$ be given by \beq \lab{eq:p}  p(g_\la, g_\mu, g_\nu, h) := Z^{-1}(a, b, c) \exp\left(-\frac{1}{2}\left(|\la|^2/a^2 + |\mu|^2/b^2 + |\nu|^2/c^2\right)\right), \eeq for some normalizing constant $Z^{-1}(a, b, c)$ independent of $\la, \mu$ and $\nu$. 
\end{definition}
Observe that the marginal on $(\la, \mu, \nu)$ is the density $\rho_{GI}(\la, \mu, \nu).$ The volume of $\mathbb{A}^3(\la, \mu; \nu)$ equals $V_n(\la)V_n(\mu)V_n(\nu)|H_n(\la, \mu; \nu)|/V_n^3(\tau_n).$
Therefore, $$Z^{-1}(a, b, c) = \left((2\pi)^{n} \left(\frac{a^2b^2c^2}{a^2 + b^2 + c^2}\right)^{n^2/2}\right)^{-1} V_n(\tau_n)^2.$$ Recall that  the differential entropy $\mathrm{ent}(q)$  of a density $q$ supported on $\mathbb{A}^3$   is given by $\int_{\mathbb{A}^3} q \ln \frac{1}{q} \Leb(dx), $ where $\Leb$ denotes the Lebesgue measure  on $\mathbb{A}^3.$
Further,  $$\E_p |\la|^2 = \frac{a^2(b^2 + c^2)n^2}{(a^2 + b^2 + c^2)}, \quad \E_p |\mu|^2 = \frac{b^2(c^2 + a^2)n^2}{(a^2 + b^2 + c^2)}.$$  and
$$\E_p |\nu|^2 = \frac{c^2(a^2 + b^2)n^2}{(a^2 + b^2 + c^2)},$$  and because $\mathrm{ent}(p) = \int_{\mathbb{A}^3} p \ln \frac{1}{p} \Leb(dx), $ which by (\ref{eq:p}), equals $$\E_p (\ln Z(a, b, c) + \frac{|\la|^2}{2a^2} + \frac{|\mu|^2}{2b^2} + \frac{ |\nu|^2}{2c^2}),$$ the differential entropy of $p$ is given by \beq \lab{entp}\mathrm{ent}(p) = n^2 + \ln Z(a, b, c).\eeq
We will use the following elementary maximum entropy principle. Unlike the standard positive-temperature Gaussian form, it only requires that the displayed exponential density be normalizable.
\begin{lemma}\lab{lem:Boltz}
Let $\mathcal{D}$ be an open set in $\R^m$, let $T_1,\ldots,T_k$ be measurable functions on $\mathcal{D}$, and let $\theta_1,\ldots,\theta_k \in \R$. Suppose
\[
Z=\int_{\mathcal{D}}\exp\left(-\sum_{j=1}^k \theta_j T_j(x)\right)dx<\infty,
\]
and let $\p$ be the probability measure on $\mathcal{D}$ with density
\[
g(x)=Z^{-1}\exp\left(-\sum_{j=1}^k \theta_j T_j(x)\right).
\]
Then among all probability distributions $\p'$ on $\mathcal{D}$ such that $\E_{\p'}T_j=\E_\p T_j$ for all $j$, the unique probability distribution that maximizes differential entropy is $\p$.
\end{lemma} 
\begin{proof}
It is enough to consider competitors $\p'$ that are absolutely continuous with respect to Lebesgue measure and have finite entropy; singular competitors do not improve the supremum. Let $f$ be the density of such a competitor, and write $h(f)=\int_{\mathcal{D}} f(x)\log(1/f(x))dx$.

Consider the Kullback–Leibler divergence between the two distributions:

\[
D_{KL}(f \parallel g) = \int_{\mathcal{D}} f(x) \log \left( \frac{f(x)}{g(x)} \right) \, dx \geq 0,
\]
with equality if and only if $f=g$ almost everywhere. Since the expectations of the sufficient statistics agree,
\[
\int_{\mathcal{D}} f(x)\log g(x)\,dx
=-\log Z-\sum_{j=1}^k\theta_j\E_{\p'}T_j
=-\log Z-\sum_{j=1}^k\theta_j\E_{\p}T_j
=\int_{\mathcal{D}} g(x)\log g(x)\,dx.
\]
Thus \(D_{KL}(f\parallel g)=h(g)-h(f)\), proving the claim.
\end{proof}
As an immediate consequence, we have the following.

\begin{thm}\lab{thm:gaussian1}
    For $\la, \mu, \nu \in \Spec_n + \R\one$, (where $\one$ denotes the vector of all ones),
among all densities $q$ supported on $\mathbb{A}^3$ such that $$\E_q |\la|^2 = \frac{a^2(b^2 + c^2)n^2}{(a^2 + b^2 + c^2)}, \quad \E_q |\mu|^2 = \frac{b^2(c^2 + a^2)n^2}{(a^2 + b^2 + c^2)}.$$  and
$$\E_q |\nu|^2 = \frac{c^2(a^2 + b^2)n^2}{(a^2 + b^2 + c^2)},$$ 
the probability density $p$ given by (\ref{eq:p}) is the unique probability density that maximizes the differential entropy. 
(When  $- c^2 > a^2 + b^2$, but, $a, b > 0$, the same statement holds.)
\end{thm}
\begin{proof}
Apply Lemma~\ref{lem:Boltz} on $\mathbb{A}^3$ with sufficient statistics $|\la|^2$, $|\mu|^2$, and $|\nu|^2$, and coefficients $1/(2a^2)$, $1/(2b^2)$, and $1/(2c^2)$. In the imaginary case the third coefficient is negative, but the partition function is finite by Theorem~\ref{thm:gaussian}; the moment identities are those computed above and in Lemma~\ref{lem:1}.
\end{proof}

 \section{Bounds on the volume  of a hive polytope}\lab{sec:lim-log-vol}

Recall that we set
   \[
   \Gamma(t) = \frac{2}{\pi} \int_{-1}^{t} \sqrt{1 - x^2} \, dx, \quad -1 \leq t \leq 1,
   \]
   and \( t = t(k, n) = \Gamma^{-1}(k/n) \) where \( k = k(n) \) is such that \( k/n \to a \in (0,1) \) as \( n \to \infty \). If \( x_k \) denotes the \( k \)-th eigenvalue in the GUE, then by a theorem of Gustavsson (see \cite[Theorem 1.1]{Gustavsson}), it holds that as \( n \to \infty \)
   \[
   \frac{x_k/\sqrt{2} - t \sqrt{2 n}}{\sqrt{\frac{\ln n}{4(1-t^2)n}}} \xrightarrow{d} N(0,1)
   \]
   in distribution. We have $x_k/\sqrt{2}$ in the above expression rather than $x_k$ as in \cite[Theorem 1.1]{Gustavsson}, because according to our convention, the standard deviation of the real and imaginary parts of a non-diagonal entry of a GUE matrix is $\frac{1}{\sqrt{2}}$ rather than $\frac{1}{2}$.

\begin{lemma}[Eigenvalue rigidity]\lab{lem:rigidity}
Let \(A\) have the distribution of GUE, in the normalization used in this paper.
Let
\[
\lambda_1 \leq \lambda_2 \leq \cdots \leq \lambda_n
\]
denote its eigenvalues, and let \(\gamma_i\) be the corresponding classical location for the semicircle law, defined by
\[
\int_{-\infty}^{\gamma_i}
\frac{1}{2\pi}\sqrt{(4-x^2)_+}\,dx
=
\frac{i}{n}.
\]
Then there are absolute constants \(c,C>0\) such that for every \(T>0\),
\[
\mathbb P\left(
n^{1/6}\min(i,n-i+1)^{1/3}
\left|\lambda_i-\sqrt n\,\gamma_i\right|
\geq T
\right)
\leq
C n^C\exp(-T^c)
\]
uniformly in \(i\).  In particular, for a suitable fixed power \(L\), with probability \(1-o(1)\),
\[
\left|\lambda_i-\sqrt n\,\gamma_i\right|
\leq
n^{-1/6}\min(i,n-i+1)^{-1/3}(\log n)^L
\]
for all \(1\leq i\leq n\).
\end{lemma}
\begin{proof}
This is the Tao--Vu eigenvalue rigidity estimate \cite{TaoVu}, specialized to the GUE normalization used here.
\end{proof}
   
\noindent Recall {\bf Definition}~\ref{def:10}:\\       Let $\la^{\mathrm{cl}}:[0, 1] \ra \R$ be given by
   $\la^{\mathrm{cl}}(x) = 2\Gamma^{-1}(1 - x)$ for $x \in [0, 1]$. Note that the discrete counterpart as specified in Definition~\ref{def:lan} satisfies  $\la_n^{\mathrm{cl}}(i) = 2n^2 \int_{1 - \frac{i}{n}}^{1 - \frac{i-1}{n}}\Gamma^{-1}(x)dx$ for $i \in [n]$.
    
   \begin{definition}\lab{def:11}
     Let $\tau:[0, 1] \ra \R$ be given by
    $\tau_n(x) = \frac{1}{2}-x$ for $x \in [0, 1]$. Recall that $\tau_n := (\tau_n(1), \dots, \tau_n(n))$ denotes the vector whose $i^{th}$ coordinate  is given by $\tau_n(i) =  \frac{n + 1}{2} - i$ for $i \in [n]$, which is in accordance with Definition~\ref{def:lan}.
    \end{definition}

   We are interested in the logarithm of the volume of the hive polytope whose boundary eigenvalues are $(\sigma_\la \lambda_n^{\mathrm{cl}}, \sigma_\mu \lambda_n^{\mathrm{cl}},  \sigma_\nu \lambda_n^{\mathrm{cl}})$. 
   \begin{lemma}\lab{lem:gamma} Let $\la$ be any bounded, monotonically decreasing function defined on $[0, 1]$.
   The limit \beq\lab{eq:equality1} \ln V(\lambda) := \lim_{n \ra \infty} \frac{2}{n^2} \ln \left(\frac{V_n(\lambda_n)}{V_n(\tau_n)}\right) \eeq exists and when $\la = \la^{\mathrm{cl}}$, we have
  \beq \lab{eq:equality2} \ln V(\lambda^{\mathrm{cl}}) = \frac{1}{ 4 \pi^2}\int_{-2}^{2} \int_{-2}^{2} \sqrt{4 - x^2} \sqrt{4 - y^2} \ln |x - y| dx dy - \int_0^1 \int_0^1 \ln |x - y| dx dy = \frac{5}{4}.\eeq
   \end{lemma} 

   \begin{proof} The first equality (\ref{eq:equality1}) is a result in \cite{NarSheff} (see Lemma 23 in that paper).
   Noting that $$\ln \left(\frac{V_n(\lambda_n^{\mathrm{cl}})}{V_n(\tau_n)}\right) = \left(\sum_{1 \leq i < j \leq n} \ln (\lambda_n^{\mathrm{cl}}(i) - \lambda_n^{\mathrm{cl}}(j))\right)  - \sum_{1 \leq i < j \leq n}\ln (j - i),$$
    the identity (\ref{eq:equality2}) follows from Riemann integration and Definition~\ref{def:10}, taking some care to handle the logarithmic singularities in the double integral, as was done in Lemma 23 of \cite{NarSheff}. 

    For the explicit evaluation of \begin{equation}\label{eq:2.5} \frac{\int_{-2}^{2} \int_{-2}^{2} \sqrt{4 - x^2} \sqrt{4 - y^2} \ln |x - y| dx dy}{ 4 \pi^2} = - \frac{1}{4}\end{equation} we refer the reader to \cite[Exercise 2.6.4]{AndersonGuionnetZeitouni}. The second integral \begin{equation}\label{eq:second} \int_0^1 \int_0^1 \ln |x - y| dx dy = -\frac{3}{2}\end{equation} follows easily once we note that it equals $2\int_0^1 \int_0^y \ln (y-x) dy dx$. 
\end{proof}

In the remainder of this paper, we will assume that $\sigma_\la, \sigma_\mu$ and $\sigma_\nu$ are fixed and in the acute or right angled case  $a, b, c,$ satisfy  $$\sigma_\la^2 :=  \frac{a^2(b^2 + c^2)}{(a^2 + b^2 + c^2)}, \quad \sigma_\mu^2:= \frac{b^2(c^2 + a^2)}{(a^2 + b^2 + c^2)},$$  and
$$\sigma_\nu^2 := \frac{c^2(a^2 + b^2)}{(a^2 + b^2 + c^2)}.$$  
If $\sigma_\nu^2 > \sigma_\la^2 + \sigma_\mu^2,$ then we are in the obtuse case, and $a, b > 0$ and $-c^2 > a^2 + b^2.$ However the same formulae hold in this case as well.
\begin{thm}\lab{thm:gaussian2}
Suppose that $\sigma_\la, \sigma_\mu$ and $\sigma_\nu$ are fixed and form a triangle with positive area. 
 Suppose $|\la| = \sigma_\la n, |\mu| = \sigma_\mu n$ and $|\nu| = \sigma_\nu n,$ where $\la, \mu, \nu \in \Spec_n + \R \one$. Then,
$$\frac{2}{n^2} \ln \left(\frac{V_n(\sqrt{n}\lambda)V_n(\sqrt{n}\mu)V_n(\sqrt{n}\nu)}{V_n(\tau_n)^3} |H_n(\sqrt{n}\lambda, \sqrt{n}\mu; \sqrt{n}\nu)|\right) \leq $$
$$5 +  \left( \ln \frac{a^2b^2c^2}{a^2 + b^2 + c^2}\right) + o(1).$$ 
\end{thm}
\begin{proof}
Without loss of generality assume that all angles are acute, except possibly the angle opposite the side of length $\sigma_\nu.$  
Let $\mathbb{T}$ be the set of all $\tilde{\la}, \tilde{\mu}, \tilde{\nu}$ such that $\tilde{\la} -\la, \tilde{\mu} - \mu, \tilde{\nu} - \nu \in \Spec_n + \R \one$ and $H_n(\tilde{\la} - \la, \tilde{\mu} - \mu; \tilde{\nu} - \nu) \neq \emptyset$ and $|\tilde{\la} - \la|^2 + |\tilde{\mu} - \mu|^2 + |\tilde{\nu} - \nu|^2 \leq n^{-1}.$ Consider the restriction of $\rho_{GI}$ to $\mathbb{T}$ normalized to be a probability density, and let $\rho_{GI}^{\mathbb{T}}$ denote this restriction.
For any $\tla, \tmu, \tnu \in \mathbb{T},$ we have $$V(\tla)V(\tmu)V(\tnu)|H_n(\tla, \tmu; \tnu)| \geq V(\la)V(\mu)V(\nu)|H_n(\la, \mu; \nu)|.$$
Secondly $(\vol_{3n-1}\mathbb{T})^{\frac{1}{n^2}} = 1 +o_n(1),$ because the $3n-1$ dimensional convex set $\mathbb{T}$ is contained in a ball of radius $n^{-\frac{1}{2}}$ centered at $(\la, \mu, \nu),$ and contains a ball of radius $n^{-C}$ for some positive constant $C$. This latter fact can be proved by considering as the center of the ball, a hive in which all the rhombus inequalities are satisfied with a slack of at least $n^{-C + 10}$. For an example of such a hive, take the values of the hive at different positions to be given by a quadratic function of two spatial variables. Let $q$ be the measure on $\phi^{-1}(\rho_{GI}^{\mathbb{T}})$ whose density at a point $x$ in its support is proportional to $p$ at $x.$ Then, the differential entropy of $q$  $$ \mathrm{ent}(q) \geq \ln\left(\frac{V_n(\la)V_n(\mu)V_n(\nu)}{V_n(\tau_n)^3} |H_n(\la, \mu; \nu)|\right),$$ and also $$\E_q |\la|^2 \leq  \sigma_\la^2 n (1 + o(1)), \quad \E_q |\mu|^2 = \sigma_\mu^2 n (1 + o(1)), \quad \E_q |\nu|^2 = \sigma_\nu^2 n (1 + o(1)).$$
Applying Theorem~\ref{thm:gaussian1} to $q$ we see that the differential entropy of $q$ is at most $o(n^2)$ plus the differential entropy of $p$ and so, 
\beqs  \ln\left(\frac{V_n(\la)V_n(\mu)V_n(\nu)}{V_n(\tau_n)^3} |H_n(\la, \mu; \nu)|\right) \leq \mathrm{ent}(q) & \leq & \\
\ent(p) = n^2 + \frac{n^2}{2}\ln\left(\frac{a^2b^2c^2}{a^2 + b^2 + c^2}\right) + o(n^2) - \ln (V_n(\tau_n)^2) & \leq & \\
n^2 + \frac{n^2}{2}\ln\left(\frac{a^2b^2c^2}{a^2 + b^2 + c^2}\right) + o(n^2) - n^2(- \frac{3}{2} + \ln n).\eeqs
Therefore, $$\frac{2}{n^2}\ln\left(\frac{V_n(\sqrt{n}\la)V_n(\sqrt{n}\mu)V_n(\sqrt{n}\nu)}{V_n(\tau_n)^3} |H_n(\sqrt{n}\la, \sqrt{n}\mu; \sqrt{n}\nu)|\right) \leq $$
 $$5 +  \left( \ln \frac{a^2b^2c^2}{a^2 + b^2 + c^2}\right) + o(1).$$

\end{proof}

\begin{thm}[Expected GUE boundary in the positive raw regime]\lab{thm:gaussian3-positive}
Assume that the associated raw parameters satisfy \(a,b,c>0\). Let
\(\overline{\la}_{GUE,n}\) denote the vector whose \(i^{th}\) coordinate is the expected value of the \(i^{th}\) eigenvalue of a GUE matrix, listed in non-increasing order.  Set
\[
\lambda_n=\sqrt n\,\sigma_\la \overline{\la}_{GUE,n},\qquad
\mu_n=\sqrt n\,\sigma_\mu \overline{\la}_{GUE,n},\qquad
\nu_n=\sqrt n\,\sigma_\nu \overline{\la}_{GUE,n}.
\]
Then
\[
\frac{2}{n^2}\log\left(
\frac{V_n(\lambda_n)V_n(\mu_n)V_n(\nu_n)}{V_n(\tau_n)^3}
|H_n(\lambda_n,\mu_n;\nu_n)|
\right)
=
5+\log\left(\frac{a^2b^2c^2}{a^2+b^2+c^2}\right)+o(1).
\]
\end{thm}
\begin{proof}
In the positive raw-parameter regime the quadratic Gaussian factor in \(F\) is log-concave.  Together with the log-concavity of \(V_n\) on the Weyl chamber and the log-concavity of hive volumes, this makes \(F\) log-concave on its support.  The normalized density \(F/Z_n\) has mean
\[
\left(\frac{\lambda_n}{\sqrt n},\frac{\mu_n}{\sqrt n},\frac{\nu_n}{\sqrt n}\right).
\]
Fradelizi's theorem therefore gives the value of \(F\) at this mean within an \(e^{-O(n)}\) factor of the supremum of \(F\).  Since \(e^{-O(n)}=e^{-o(n^2)}\), the normalization calculation is the same as in the proof of Theorem~\ref{thm:gaussian3}.
\end{proof}

\begin{remark}
Theorem~\ref{thm:gaussian3-positive} is not used in the obtuse case.  Its proof is included to record the original Fradelizi argument in exactly the regime where \(F\) is log-concave.
\end{remark}

\begin{thm}[Rigidity-box lower bound]\lab{thm:gaussian3}
Let \(L\) be a sufficiently large fixed constant, and put
\[
r_{i,n}=n^{1/3}\min(i,n-i+1)^{-1/3}(\log n)^L.
\]
Then there exist boundary triples
\[
(\Lambda_n^*,M_n^*,N_n^*)\in(\Spec_n+\mathbb R\one)^3,
\qquad
\sum_i\Lambda_n^*(i)+\sum_iM_n^*(i)=\sum_iN_n^*(i),
\]
such that, in each of the three spectra, the minimum spacing between consecutive
eigenvalues is at least \(n^{-100}\), the hive polytope
\(H_n(\Lambda_n^*,M_n^*;N_n^*)\) contains a Euclidean ball of radius
\(n^{-100}\) in its affine span, and for all \(i\),
\[
\left|\Lambda_n^*(i)-\sigma_\la\lambda_n^{\mathrm{cl}}(i)\right|
\leq C r_{i,n},\quad
\left|M_n^*(i)-\sigma_\mu\lambda_n^{\mathrm{cl}}(i)\right|
\leq C r_{i,n},
\]
and
\[
\left|N_n^*(i)-\sigma_\nu\lambda_n^{\mathrm{cl}}(i)\right|
\leq C r_{i,n},
\]
where \(C\) depends only on \(\sigma_\la,\sigma_\mu,\sigma_\nu\), and such that
\[
\frac{2}{n^2} \ln \left(\frac{V_n(\Lambda_n^*)V_n(M_n^*)V_n(N_n^*)}{V_n(\tau_n)^3}
|H_n(\Lambda_n^*,M_n^*;N_n^*)|\right)
=
5+\ln\left(\frac{a^2b^2c^2}{a^2+b^2+c^2}\right)+o(1).
\]
This holds in both raw-parameter regimes \(a,b,c>0\) and \(a,b>0,\ -c^2>a^2+b^2\).
\end{thm}

\begin{proof}
For $\lambda,\mu,\nu \in \Spec_n + \R\one$, recall
\[
F(\lambda,\mu,\nu)
=
\frac{V_n(\lambda)V_n(\mu)V_n(\nu)}{V_n(\tau_n)}
|H_n(\lambda,\mu;\nu)|
\exp\left[-\frac12\left(
\frac{|\lambda|^2}{a^2}
+\frac{|\mu|^2}{b^2}
+\frac{|\nu|^2}{c^2}\right)\right],
\]
and otherwise \(F=0\).  By Theorem~\ref{thm:gaussian},
\[
Z_n:=\int_{\sum\lambda+\sum\mu=\sum\nu}F(\lambda,\mu,\nu)\,d\lambda\,d\mu\,d\nu
=
(2\pi)^n
\left(\frac{a^2b^2c^2}{a^2+b^2+c^2}\right)^{n^2/2}.
\]
Let \(\rho_{GI}=Z_n^{-1}F\).  The one-boundary marginals of \(\rho_{GI}\) are the scaled GUE spectral measures with variances \(\sigma_\la^2,\sigma_\mu^2,\sigma_\nu^2\), respectively; this is the computation preceding Lemma~\ref{lem:1}, together with Lemma~\ref{lem:1} in the imaginary-parameter case.

Let \(\mathcal B_n\) be the rigidity box displayed in the statement, at hive-boundary scale, and let
\[
\widetilde{\mathcal B}_n
=
\{(\lambda,\mu,\nu):(\sqrt n\,\lambda,\sqrt n\,\mu,\sqrt n\,\nu)\in\mathcal B_n\}.
\]
Thus \(\widetilde{\mathcal B}_n\) is centered at
\[
n^{-1/2}(\sigma_\la\lambda_n^{\mathrm{cl}},
\sigma_\mu\lambda_n^{\mathrm{cl}},
\sigma_\nu\lambda_n^{\mathrm{cl}})
\]
and has coordinate widths \(n^{-1/2}r_{i,n}=n^{-1/6}\min(i,n-i+1)^{-1/3}(\log n)^L\).  By Lemma~\ref{lem:rigidity} and a union bound over the three GUE marginals,
\[
\int_{\widetilde{\mathcal B}_n}F\geq (1-o(1))Z_n.
\]
Moreover,
\[
\log\operatorname{vol}_{3n-1}(\widetilde{\mathcal B}_n)=O(n\log n)=o(n^2),
\]
because \(\sum_i\log r_{i,n}=O(n\log n)\).  Choose
\((\lambda_n^*,\mu_n^*,\nu_n^*)\in\widetilde{\mathcal B}_n\) at which \(F\) is maximal.  Then
\[
\log F(\lambda_n^*,\mu_n^*,\nu_n^*)
\geq \log Z_n-o(n^2).
\]
The rigidity bounds imply
\[
|\lambda_n^*|^2=\sigma_\la^2n^2+o(n^2),\quad
|\mu_n^*|^2=\sigma_\mu^2n^2+o(n^2),\quad
|\nu_n^*|^2=\sigma_\nu^2n^2+o(n^2),
\]
and therefore
\[
\frac12\left(
\frac{|\lambda_n^*|^2}{a^2}
+\frac{|\mu_n^*|^2}{b^2}
+\frac{|\nu_n^*|^2}{c^2}\right)
=n^2+o(n^2).
\]
It follows that
\[
\frac{2}{n^2}\log\left(
\frac{V_n(\lambda_n^*)V_n(\mu_n^*)V_n(\nu_n^*)}{V_n(\tau_n)}
|H_n(\lambda_n^*,\mu_n^*;\nu_n^*)|\right)
=
2+\log\left(\frac{a^2b^2c^2}{a^2+b^2+c^2}\right)+o(1).
\]
Set
\[
\Lambda_n^*=\sqrt n\,\lambda_n^*,\qquad
M_n^*=\sqrt n\,\mu_n^*,\qquad
N_n^*=\sqrt n\,\nu_n^*.
\]
Since \((\lambda_n^*,\mu_n^*,\nu_n^*)\in\widetilde{\mathcal B}_n\), this triple lies in the hive-scale box \(\mathcal B_n\), and hence satisfies the rigidity estimates in the statement.  Rescaling all three boundaries by \(\sqrt n\) multiplies the displayed product by
\[
(\sqrt n)^{3\binom n2+\binom{n-1}{2}},
\]
and using \(V_n(\tau_n)^2=\exp(-3n^2/2)n^{n^2}(1+o(1))^{n^2}\), obtained in \((\ref{eq:second})\), gives the claimed normalized formula.  The matching upper bound is Theorem~\ref{thm:gaussian2}, applied to these rigidity-box boundaries; their Euclidean norms are \(\sigma_\la n+o(n)\), \(\sigma_\mu n+o(n)\), and \(\sigma_\nu n+o(n)\), so the same proof gives the same upper bound with \(o(1)\) error.
Finally, we enforce the eigengap condition by adding a $n^{-100} \tau$ to each spectrum, and noting that the quadratic hive with each boundary equal to $\tau$ contains a constant sized  Euclidean ball in it.
\end{proof}

  \begin{lemma} \lab{lem:gaussian4} Let $ \overline{\la}_{GUE, n}$ denote the vector whose $i^{th}$ coordinate is the expected value of the $i^{th}$ eigenvalue of a GUE matrix with diagonal entries that are $N(0, 1)_\R$ and non-diagonal entries that are $N(0, 1)_\mathbb{C}$, listed in non-increasing order. 
   Then, $$V_n\left(\overline{\la}_{GUE, n}\right)^\frac{2}{n^2} = \exp\left(-\frac{1}{4} - o(1)\right)\sqrt{n}.$$
    \end{lemma}
    Note: here $o(1)$ denotes a term whose {\it magnitude} tends to zero as $n$ tends to infinity.
\begin{proof} Let $\gamma$ be the uniform measure on $B := \{x \in \Spec_n + \R\one| |x| \leq n^4\}$ of radius $n^4$ in $\Spec_n + \R \one$, centered at the origin whose density is $Z^{-1} I(B)$.
Recall from the proof of Theorem~\ref{thm:gaussian} that the density of the GUE eigenvalues is given by 
$$GUE(\lambda, 1) := \left(\frac{1}{(2\pi)^\frac{n}{2}}\right)  \frac{V_n(\la)^2}{V_n(\tau_n)}\exp\left(-\frac{|\la|^2}{2}\right)$$
for $\la \in \Spec_n + \R\one$ and $GUE(\la, 1) = 0$ otherwise. 

Note that  \beqs \sup_{\Spec_n + \R \one} \left(\frac{V_n(\la)^2}{V_n(\tau_n)}\right)^{\frac{2}{n^2}}\exp\left(-\frac{|\la|^2}{n^2}\right) =  \sup_{B} \left(\frac{V_n(\la)^2}{V_n(\tau_n)}\right)^{\frac{2}{n^2}}\exp\left(-\frac{|\la|^2}{n^2}\right).\eeqs and that 
$$ \frac{1}{n!} \leq vol(B ) \leq n^{4n}.$$ Thus, 
\beqs \left(\sup_{B} \left(\frac{1}{2\pi}\right)^{\frac{n}{2}}\frac{V_n(\la)^2}{V_n(\tau_n)} \exp\left(- \frac{|\la|^2}{2}\right)\right)^\frac{2}{n^2}\eeqs is greater or equal to  \beqs \left(\E_{\gamma} \left(\frac{1}{2\pi}\right)^{\frac{n}{2}}\frac{V_n(\la)^2}{V_n(\tau_n)} \exp\left(-\frac{|\la|^2}{2}\right)\right)^\frac{2}{n^2} \eeqs which is in turn greater or equal to \beqs
 (1 - o(1))\left(\int_B \left(\frac{1}{2\pi}\right)^{\frac{n}{2}}\frac{V_n(\la)^2}{V_n(\tau_n)} \exp\left(-\frac{|\la|^2}{2}\right) d\la \right)^\frac{2}{n^2} 
  =   (1 - o(1)).\eeqs

Therefore, \beq \sup_{\Spec_n + \R \one} \left(\frac{V_n(\la)^2}{V_n(\tau_n)}\right)^{\frac{2}{n^2}}\exp\left(-\frac{|\la|^2}{n^2}\right) & \geq & 1 - o(1).\eeq

Next, proceeding by the same Fradelizi argument as in Theorem~\ref{thm:gaussian3-positive}, we see that applying the theorem of Fradelizi (see Theorem~\ref{thm:frad}) \cite{Fradelizi}, 
\beqs \left(\frac{V_n(\overline{\la}_{GUE, n})^2}{V_n(\tau_n)}\right)^{\frac{2}{n^2}}\exp\left(-\frac{|\overline{\la}_{GUE, n}|^2}{n^2}\right) & \geq & (1 - o(1))\sup_{\Spec_n + \R \one} \left(\frac{V_n(\la)^2}{V_n(\tau_n)}\right)^{\frac{2}{n^2}}\exp\left(-\frac{|\la|^2}{n^2}\right)\\
& \geq & 1 - o(1).\eeqs
Since $|\overline{\la}_{GUE, n}|=n(1 + o(1))$ and it has already been shown in (\ref{eq:second})  that $V_n(\tau_n)^\frac{2}{n^2} = (\exp(-3/2))n(1 + o(1))$,

we see that $$V_n(\overline{\la}_{GUE, n})^\frac{2}{n^2} \geq \exp\left(-\frac{1}{4} - o(1)\right)\sqrt{n}.$$

{\color{black}{

Next, note that by (\ref{eq:2.5}), $$V_n(n^{-\frac{1}{2}}\la^{\mathrm{cl}}_n)^\frac{2}{n^2} = \exp\left(-\frac{1}{4} - o(1)\right)\sqrt{n}.$$
The minimum spacing between consecutive eigenvalues comprising the vector $\la^{\mathrm{cl}}_n$ is at least $c$. Therefore, there is a ball $\tilde{B}_n$ of radius greater than  $c\eps_0  n^{- \frac{1}{2}}$ centered at $\lacl_n$, such that for all $\la$ in this ball, 
$$V_n(n^{-\frac{1}{2}}\la)^\frac{2}{n^2} > \exp\left(-\frac{1}{4} - o(1)\right)\sqrt{n}(1 - C \eps_0).$$

For the sake of contradiction, suppose now that for some positive $\eps_1$, there is an increasing infinite sequence $n_1, n_2, \dots$ such that for all positive $i$, 
$$V_{n_i}(\overline{\la}_{GUE, n_i})^\frac{2}{n_i^2} \geq \exp\left(-\frac{1}{4} + \eps_1\right)\sqrt{n_i}.$$ 

Let $K_i \subseteq \R^{n_i}$ be the convex hull of $\{\overline{\la}_{GUE, n_i}\}\cup \tilde{B}_{n_i}.$
Then, there exists a sufficiently small universal constant $\eps_2$ such that for all $\la_{n_i}$ in the set $(\overline{\la}_{GUE, n_i})\eps_2 + (1 - \eps_2) \tilde{B}_{n_i},$ it is true by log-concavity that 
$$V_{n_i}(\la_{n_i})^\frac{2}{n_i^2} \exp(-\frac{\la_{n_i}^2}{2n_i^2})\geq \exp\left(-\frac{3}{4} + \eps_2\right)\sqrt{n_i}.$$ 
But the radius of $\tilde{B}_{n_i}$ is large enough that this by an application of (\ref{eq:second}) implies that for all sufficiently large $i$, 
$$\int_{(\overline{\la}_{GUE, n_i})\eps_2 + (1 - \eps_2) \tilde{B}_{n_i}} \left(\frac{1}{2\pi}\right)^{\frac{n}{2}}\frac{V_{n_i}(\la)^2}{V_{n_i}(\tau_{n_i})} \exp\left(-\frac{|\la|^2}{2}\right) d\la > 1,$$ which is a contradiction.
}}
\end{proof}

\noindent{\bf Ordered Mehta integral.}
We will use the following normalization of the ordered Mehta integral, equivalently the normalization of the ordered GUE eigenvalue density in \cite[Equation 2.4.23]{AndersonGuionnetZeitouni}.  For every \(s>0\),
\begin{equation}\lab{eq:ordered-mehta}
\int_{\Spec_n+\mathbb R\one}
V_n(x)^2\exp\left(-\frac{|x|^2}{2s^2}\right)\,dx
=(2\pi)^{n/2}s^{n^2}V_n(\tau_n).
\end{equation}
Equivalently, the ordered GUE eigenvalue density with variance parameter \(s^2\) is
\[
\frac{1}{(2\pi)^{n/2}s^{n^2}V_n(\tau_n)}
V_n(x)^2\exp\left(-\frac{|x|^2}{2s^2}\right),
\qquad x\in\Spec_n+\mathbb R\one.
\]
In particular, applying Fradelizi's theorem to this ordered density gives the Vandermonde upper bound
\begin{equation}\lab{eq:mehta-vandermonde-bound}
\frac{2}{n^2}\log\frac{V_n(x_n)}{V_n(\tau_n)}
\leq \frac54+\log\sigma+o(1)
\end{equation}
whenever \(x_n\in\Spec_n+\mathbb R\one\) and \(|x_n|=\sigma n^{3/2}+o(n^{3/2})\).
Here is the calculation.  Let
\[
G_s(x)=V_n(x)^2\exp\left(-\frac{|x|^2}{2s^2}\right).
\]
The normalized version of \(G_s\) in \((\ref{eq:ordered-mehta})\) is log-concave, and its mean is
\[
\overline x_s=s\,\overline{\lambda}_{GUE,n},
\]
where \(\overline{\lambda}_{GUE,n}\) is the vector of expected ordered GUE eigenvalues from Lemma~\ref{lem:gaussian4}.  Fradelizi's theorem, applied in dimension \(n\), gives
\[
\sup_{\Spec_n+\mathbb R\one}G_s(x)\leq e^nG_s(\overline x_s).
\]
Thus
\[
V_n(x_n)^2
\leq
e^nV_n(\overline x_s)^2
\exp\left(\frac{|x_n|^2-|\overline x_s|^2}{2s^2}\right).
\]
Taking \(s=\sigma\sqrt n\), dividing by \(V_n(\tau_n)^2\), and using Lemma~\ref{lem:gaussian4} together with the homogeneity of the Vandermonde determinant gives
\[
\frac{2}{n^2}\log\frac{V_n(\overline x_s)}{V_n(\tau_n)}
=\frac54+\log\sigma+o(1).
\]
The accompanying norm estimate \(|\overline{\lambda}_{GUE,n}|=n(1+o(1))\) and the hypothesis on \(x_n\) also give
\[
\frac{|x_n|^2-|\overline x_s|^2}{2s^2n^2}=o(1),
\]
so \((\ref{eq:mehta-vandermonde-bound})\) follows.


Recall that  $\De(\sigma_\la, \sigma_\mu, \sigma_\nu)$ denotes the area of the Euclidean triangle with sides $\sigma_\la, \sigma_\mu, \sigma_\nu$.  Assuming that such a triangle exists, by Heron's formula, we have 
\begin{eqnarray}\lab{eq:heron} \De(\sigma_\la, \sigma_\mu, \sigma_\nu) = \sqrt{\mathbf{s}(\mathbf{s} - \sigma_\la)(\mathbf{s} - \sigma_\mu)(\mathbf{s} - \sigma_\nu)},\end{eqnarray}
where $\mathbf{s} = (\sigma_\la + \sigma_\mu + \sigma_\nu)/2$.

\begin{lemma}\lab{lem:3}
Let $a, b, c > 0$, or $\sqrt{-1} c \in \R$ with $a, b > 0$ and $- c^2 > a^2 + b^2$. Recall that $$\sigma_\la^2 :=  \frac{a^2(b^2 + c^2)}{(a^2 + b^2 + c^2)}, \quad \sigma_\mu^2:= \frac{b^2(c^2 + a^2)}{(a^2 + b^2 + c^2)},$$  and
$$\sigma_\nu^2 := \frac{c^2(a^2 + b^2)}{(a^2 + b^2 + c^2)}.$$  Then,
$$\frac{a^2b^2c^2}{a^2 + b^2 + c^2} = 4\De(\sigma_\la, \sigma_\mu, \sigma_\nu)^2.$$
\end{lemma}
\begin{proof}
We shall expand $16 \De(\sigma_\la, \sigma_\mu, \sigma_\nu)^2$ using (\ref{eq:heron}) and show that it equals $4a^2b^2c^2/(a^2 + b^2 + c^2).$ Indeed, we have
$$16 \De(\sigma_\la, \sigma_\mu, \sigma_\nu)^2 = $$ $$(\sigma_\la + \sigma_\mu + \sigma_\nu)(\sigma_\la + \sigma_\mu - \sigma_\nu)(\sigma_\la - \sigma_\mu + \sigma_\nu)(-\sigma_\la + \sigma_\mu + \sigma_\nu).$$
This simplifies to $$((\sigma_\la + \sigma_\mu)^2 - \sigma_\nu^2)(\sigma_\nu^2 - (\sigma_\la - \sigma_\mu)^2) = $$
$$\left(\frac{2a^2b^2}{a^2 + b^2 + c^2} + \frac{2ab\sqrt{b^2 + c^2}\sqrt{c^2 + a^2}}{a^2 + b^2 + c^2}\right)\left( - \frac{2a^2b^2}{a^2 + b^2 + c^2} + \frac{2ab\sqrt{b^2 + c^2}\sqrt{c^2 + a^2}}{a^2 + b^2 + c^2}\right) = $$
$$\frac{4a^2b^2}{(a^2 + b^2 + c^2)^2}\left(c^4 + a^2 c^2 + b^2 c^2\right) = 4a^2b^2c^2/(a^2 + b^2 + c^2).$$
\end{proof}
\begin{thm}\lab{thm:5}
    Let $\la_n = \sigma_\la \lambda_n^{\mathrm{cl}}$ and $\mu_n =  \sigma_\mu \lambda_n^{\mathrm{cl}}$. Let $\g$ satisfy $\partial^-\gamma = \nu$ and for any natural number $n$ and $1 \leq i \leq n$, $\nu_n(i) = \sigma_\nu\la_n^{\mathrm{cl}}(i)$.  Let $X_n^\text{cl} = \sigma_\la U \diag(\lambda_n^{\mathrm{cl}}) U^*$ and $Y_n^\text{cl} = \sigma_\mu V \diag(\lambda_n^{\mathrm{cl}}) V^*$ where $U$ and $V$ are independent matrices from the Haar measure on the unitary group $\mathbb{U}_n$. Let $\sigma_\la, \sigma_\mu, \sigma_\nu \in \R$, such that there exists a Euclidean triangle with these sides and  $ \De  (\sigma_\la, \sigma_\mu, \sigma_\nu) > 0$. Then,
$$\lim_{\eps \ra 0} \lim_{n \ra \infty}  (2/n^2) \ln\left(\int\limits_{B_\I^n(\nu_n, \eps n^2)} \left(\frac{V_n(\nu'_n)}{V_n(\nu_n)}\right) |H_n(\la_n, \mu_n; \nu'_n)|\Leb_{n-1, 0}(d\nu'_n)\right) =$$ $$ \frac{5}{4} + \ln \left(\frac{4 \De  (\sigma_\la, \sigma_\mu, \sigma_\nu)^2}{\sigma_\la \sigma_\mu \sigma_\nu}\right).$$  \end{thm}
\begin{proof}

Recall that 
the volume of the $n^2 - n$ dimensional polytope 
\begin{eqnarray*} |G_n(\la_n, \mu_n; \nu_n, \eps n^2) | & = & \int\limits_{B_\I^n(\nu_n, \eps n^2)}|A_n(\la_n, \mu_n; \nu'_n)|\Leb_{n-1, 0}(d\nu'_n).\lab{eq:G}\\
& = & \int\limits_{B_\I^n(\nu_n, \eps n^2)} \left(\frac{V_n(\nu'_n)}{V_n(\tau_n)}\right) |H_n(\la_n, \mu_n; \nu'_n)|\Leb_{n-1, 0}(d\nu'_n).\end{eqnarray*} 
For all positive $\epsilon$, by Theorem~\ref{thm:1} we have
\begin{eqnarray*} \p_n\left[\spec(X_n^\text{cl} + Y_n^\text{cl}) \in B_\I^n(\nu_n, \eps n^2)\right] =\end{eqnarray*} \begin{eqnarray} \label{eq:3.3}\left( \frac{V_n(\tau_n)^2}{V_n(\la_n)V_n(\mu_n)} \right)\int\limits_{B_\I^n(\nu_n, \eps n^2)} \left(\frac{V_n(\nu'_n)}{V_n(\tau_n)}\right) |H_n(\la_n, \mu_n; \nu'_n)|\Leb_{n-1, 0}(d\nu'_n).\lab{eq:3.6}\end{eqnarray}

Therefore, the LHS in the statement of the theorem equals
\begin{eqnarray}\lab{eq:expr1}\lim_{\eps \ra 0} \lim_{n \ra \infty} \frac{2}{n^2} \ln \left(\frac{V_n(\la_n)V_n(\mu_n)}{V_n(\tau_n) V_n(\nu_n)}  \p_n\left[\spec(X_n^\text{cl} + Y_n^\text{cl}) \in B_\I^n(\nu_n, \eps n^2)\right]\right) .\end{eqnarray}

By Lemma~\ref{lem:gamma}, (\ref{eq:expr1}) equals
\begin{eqnarray}\lab{eq:expr2} \frac{5}{4} + \ln \frac{\sigma_\la \sigma_\mu}{\sigma_\nu}  + \lim_{\eps \ra 0} \lim_{n \ra \infty} \frac{2}{n^2} \ln \left(\p_n\left[\spec(X_n^\text{cl} + Y_n^\text{cl}) \in B_\I^n(\nu_n, \eps n^2)\right]\right) .\end{eqnarray} 
The upper bound for the expression in the statement follows directly from Theorem~\ref{thm:gaussian2}, Lemma~\ref{lem:3}, and Lemma~\ref{lem:gamma}.  We prove the matching lower bound.
As in the proof of \cite[Lemma 26]{NarSheff}, 
by the Ky Fan inequalities, (see \cite[Equation 4
(1.55)]{TaoBook}), which state that, for each $k \in [n]$, for $n \times n$ Hermitian matrices $X$ and $Y$, the sum of the top $k$ eigenvalues of $X+Y$ is less or equal to the sum of the top $k$ eigenvalues of $X$ added to the sum of the top $k$ eigenvalues of $Y,$ 
we see that $\|\spec(X + Y)\|_\I \leq \|\spec(X)\|_\I + \|\spec(Y)\|_\I.$
Let \((\Lambda_n^*,M_n^*,N_n^*)\) be the rigidity-box boundary triple from Theorem~\ref{thm:gaussian3}.  Let
\[
X_n^*=U\diag(\Lambda_n^*)U^*,\qquad
Y_n^*=V\diag(M_n^*)V^*,
\]
for the same \(U,V\) as in the statement of the theorem.  Since the side lengths in the hive-scale rigidity box are \(r_{i,n}\), their total \(\|\cdot\|_\I\)-size is \(o(n^2)\).  Hence, for every fixed \(\eps_0>0\) and all sufficiently large \(n\),
\beq \|\Lambda_n^* - \la_n\|_\I + \|M_n^* - \mu_n\|_\I + \|N_n^* - \nu_n\|_\I \leq \frac{\eps_0}{10} n^2.\eeq 
Then, by the Ky Fan inequalities, we have that
\begin{eqnarray}\lab{eq:KY2} \p_n\left[\spec(X_n^\text{cl} + Y_n^\text{cl}) \in B_\I^n(\nu_n, (\eps + \eps_0) n^2)\right] \geq \p_n\left[\spec(X_n^* + Y_n^*) \in B_\I^n(N_n^*, \eps  n^2)\right].\end{eqnarray}
Similarly, 
\begin{eqnarray}\lab{eq:KY3} \p_n\left[\spec(X_n^\text{cl} + Y_n^\text{cl}) \in B_\I^n(\nu_n, (\eps - \eps_0) n^2)\right] \leq \p_n\left[\spec(X_n^* + Y_n^*) \in B_\I^n(N_n^*, \eps n^2)\right].\end{eqnarray}

We now apply the spectral formula once more, but with the two fixed boundary
vectors \(\Lambda_n^*\) and \(M_n^*\).  Thus, for every Borel set
\(E\subset \Spec_n+\mathbb R\one\),
\[
\p_n\left[\spec(X_n^*+Y_n^*)\in E\right]
=
\frac{V_n(\tau_n)^2}{V_n(\Lambda_n^*)V_n(M_n^*)}
\int_E
\frac{V_n(\nu')}{V_n(\tau_n)}
|H_n(\Lambda_n^*,M_n^*;\nu')|
\Leb_{n-1,0}(d\nu').
\]
We use this with \(E=B_\I^n(N_n^*,\eps n^2)\).  The only point requiring
comment is the passage from the single boundary \(N_n^*\) to a positive-volume
set of nearby boundaries.  This is not an application of log-concavity to the
full three-boundary density \(F\); in the non-positive raw-parameter regime
that density need not be log-concave.

The eigengap condition in Theorem~\ref{thm:gaussian3} makes the Vandermonde
part harmless.  Choose \(K>200\), and let \(E_n\) be a convex set of admissible
third boundaries contained in the coordinate box
\[
\|\nu'-N_n^*\|_\infty\leq n^{-K}
\]
and with \(\log\Leb_{n-1,0}(E_n)=o(n^2)\).  Since the consecutive gaps of
\(N_n^*\) are at least \(n^{-100}\), every \(\nu'\in E_n\) remains in the Weyl
chamber and
\[
\log V_n(\nu')\geq \log V_n(N_n^*)-o(n^2).
\]
Indeed, each gap changes by at most \(2n^{-K}\), so the total change in
\(\log V_n\) is bounded by \(O(n^2 n^{100-K})=o(n^2)\).

It remains only to choose \(E_n\) inside the Horn cone with hive fibers not too
small.  Here we use the interior-ball conclusion of Theorem~\ref{thm:gaussian3}
directly.  Let \(h_n^*\in H_n(\Lambda_n^*,M_n^*;N_n^*)\) be the center of a
Euclidean ball of radius \(n^{-100}\) contained in this hive fiber.  If
\(\|\nu'-N_n^*\|_\infty\leq n^{-K}\) with \(K>200\), then changing the third
boundary from \(N_n^*\) to \(\nu'\) changes each affected rhombus inequality by
at most \(O(n^{-K})\).  Hence the same center \(h_n^*\), with the third boundary
replaced by \(\nu'\), still has polynomial slack, and the corresponding fiber
contains a Euclidean ball of radius \(n^{-101}\).  Thus we may take
\[
E_n=\{\nu'\in\Spec_n+\mathbb R\one:
\sum_i\nu'(i)=\sum_iN_n^*(i),\ \|\nu'-N_n^*\|_\infty\leq c n^{-K}\}
\]
with \(c>0\) small.  This set is contained in the Horn cone for the fixed pair
\((\Lambda_n^*,M_n^*)\), has \(\log\Leb_{n-1,0}(E_n)=o(n^2)\), and lies in
\(B_\I^n(N_n^*,\eps n^2)\) for all large \(n\).  The ball-volume bound,
equivalently Brunn--Minkowski log-concavity of hive fibers, gives
\[
|H_n(\Lambda_n^*,M_n^*;\nu')|
\geq e^{-o(n^2)}
|H_n(\Lambda_n^*,M_n^*;N_n^*)|
\qquad(\nu'\in E_n).
\]
Combining this with the Vandermonde estimate and integrating over \(E_n\),
\[
\int_{E_n}
\frac{V_n(\nu')}{V_n(\tau_n)}
|H_n(\Lambda_n^*,M_n^*;\nu')|
\Leb_{n-1,0}(d\nu')
\geq
e^{-o(n^2)}
\frac{V_n(N_n^*)}{V_n(\tau_n)}
|H_n(\Lambda_n^*,M_n^*;N_n^*)|.
\]
Thus replacing the point \(N_n^*\) by the
\(B_\I^n(N_n^*,\eps n^2)\)-window loses only \(e^{o(n^2)}\) in the lower-bound
estimate.  Combining this with the standard GUE Vandermonde upper bound
\[
\frac{2}{n^2}\log\frac{V_n(\Lambda_n^*)}{V_n(\tau_n)}\leq \frac54+\log\sigma_\la+o(1),
\qquad
\frac{2}{n^2}\log\frac{V_n(M_n^*)}{V_n(\tau_n)}\leq \frac54+\log\sigma_\mu+o(1),
\]
which follows from (\ref{eq:mehta-vandermonde-bound}) and the fact that \(|\Lambda_n^*|=\sigma_\la n^{3/2}+o(n^{3/2})\), \(|M_n^*|=\sigma_\mu n^{3/2}+o(n^{3/2})\), we obtain
$$\liminf_{\eps \ra 0} \liminf_{n \ra \infty}  \frac{2}{n^2} \ln \p_n\left[\spec(X_n^* + Y_n^*) \in B_\I^n(N_n^*, \eps  n^2)\right] \geq \ln \left(\frac{4 \De  (\sigma_\la, \sigma_\mu, \sigma_\nu)^2}{\sigma_\la^2 \sigma_\mu^2}\right).$$
Together with (\ref{eq:KY2}) and the upper bound noted above, this proves the theorem.
\end{proof}

\section{Bounds on the surface tension}\lab{sec:bounds_surf}

Recall from (\ref{eq:I}) that for a concave $\gamma'$ such that $\nu' = \partial^- \gamma'$,  \begin{eqnarray*} I(\g') :=  \ln\left(\frac{V(\la)V(\mu)}{V(\nu')}\right) + \inf\limits_{h' \in H(\la, \mu; \nu')}  \int_T \sigma((-1)(\hess h')_{ac})\Leb_2(dx).\end{eqnarray*}\\

\noindent {\bf Theorem~\ref{thm:6} (restated).}\\  Let $\sigma_\la, \sigma_\mu, \sigma_\nu \in \R$, such that there exists a Euclidean triangle with these sidelengths and  $ \De  (\sigma_\la, \sigma_\mu, \sigma_\nu) > 0$.  Let $\la, \mu$ and $\nu$ denote bounded,  strongly decreasing  functions from $[0, 1]$ to $\R$.  Suppose that $\|\la\|_{L^2} = \sigma_\la$, $\|\mu\|_{L^2} = \sigma_\mu$ and $\|\nu\|_{L^2} = \sigma_\nu$. 

Then,
    \beq  \inf\limits_{h \in H(\la, \mu; \nu)}  \int_T \sigma((-1)(\hess h)_{ac})\Leb_2(dx) & \geq & \nonumber \\ \ln V(\la) + \ln V(\mu) + \ln V(\nu)  - 5 & - &  \ln\left(4  \De  (\sigma_\la, \sigma_\mu, \sigma_\nu)^2 \right),\eeq
 with equality if
   $\la, \mu, \nu$ are respectively $\sigma_\la \la^{\mathrm{cl}}, \sigma_\mu \la^{\mathrm{cl}}, \sigma_\nu \la^\mathrm{cl}$.
  \begin{proof}
    By Theorem~\ref{thm:1},  \begin{eqnarray*}\lab{eq:2.4new_mod} \rho_n\left[\spec(X_n + Y_n) = \nu_n\right] =  \frac{V_n(\nu_n)V_n(\tau_n)}{V_n(\la_n)V_n(\mu_n)} |H_n(\la_n, \mu_n; \nu_n)|.\end{eqnarray*}
By Theorem~\ref{thm:4}, \begin{eqnarray*}\lim\limits_{n \ra \infty}\left(\frac{-2}{n^2}\right)\ln \p_n\left[\spec(Z_n) \in {B}_\I^n(\nu_n, \eps)\right] = \inf\limits_{\partial^-\gamma'\in {B}_\I(\nu, \eps)}I(\g').\end{eqnarray*}
Recall from (\ref{eq:I}) that \begin{eqnarray*}\lab{eq:I2} I(\g') :=  \ln\left(\frac{V(\la)V(\mu)}{V(\nu')}\right) + \inf\limits_{h' \in H(\la, \mu; \nu')}  \int_T \sigma((-1)(\hess h')_{ac})\Leb_2(dx).\end{eqnarray*}

Applying Lemma 47 of \cite{NarSheff}, which states that if $\partial^-\gamma = \nu$, where $\gamma$ is Lipschitz and concave, then $$\lim\limits_{\eps \ra 0}  \inf\limits_{\partial^-\gamma'\in {B}_\I(\nu, \eps)}I(\g') = I(\g),$$ we have
 \beqs \lim_{\eps \ra 0} \lim_{n \ra \infty}  (2/n^2) \ln\left(\int\limits_{B_\I^n(\nu_n, \eps n^2)} \left(\frac{V_n(\nu'_n)V_n(\tau_n)^2}{V_n(\la_n)V_n(\mu_n)V_n(\nu_n)}\right) |H_n(\la_n, \mu_n; \nu'_n)|\Leb_{n-1, 0}(d\nu'_n)\right) = \\ - \ln V(\nu) - I(\g). \eeqs
 \begin{claim} $\limsup\limits_{\eps \ra 0} \sup\limits_{\nu'_n \in B_\I^n(\nu_n, \eps n^2)}\frac{\|\nu'_n\|}{\|\nu_n\|} = 1 $.
 \end{claim}
 \begin{proof} This is because  if a piecewise linear, concave  function on $[0, 1]$ with an a priori bounded Lipschitz norm approximates a strongly concave Lipschitz function on the same domain within $\eps$ in the sup norm,  then the $L^2$ norm of the left derivative of the approximating function is $o(1)$ as $\eps \ra 0$.
 \end{proof}
  By Theorem~\ref{thm:gaussian2} and Lemma~\ref{lem:3},
 \beq\lab{eq:@} \frac{2}{n^2} \ln \left(\frac{V_n(\lambda_n)V_n(\mu_n)V_n(\nu'_n)}{V_n(\tau_n)^3} |H_n(\lambda_n, \mu_n; \nu'_n)|\right) \leq 
 5 +  \ln \left(4 \De^2(\sigma_\la, \sigma_\mu, \sigma_\nu)\right) + o(1).\eeq
 Therefore,
  \beqs & - & \ln V(\nu) - I(\g)  \leq   -\ln(V(\la)^2 V(\mu)^2V(\nu)) \\  & + &   \lim_{\eps \ra 0} \lim_{n \ra \infty}  (2/n^2) \ln\left(\int\limits_{B_\I^n(\nu_n, \eps n^2)} \exp(\frac{n^2}{2}(5 +  \ln \left(4 \De^2(\sigma_\la, \sigma_\mu, \sigma_\nu)\right) + o(1))))\Leb_{n-1, 0}(d\nu'_n)\right).\eeqs
 Because the integrand has the form $\exp(c n^2)$ but the domain of integration is $n $ dimensional and is a convex set whose inradius and circumradius are polynomially related (as a function of  $n$), 
 we have 
 \beqs  \inf\limits_{h \in H(\la, \mu; \nu)}  \int_T \sigma((-1)(\hess h)_{ac})\Leb_2(dx)& \geq & \nonumber\\ \ln V(\la) + \ln V(\mu) + \ln V(\nu)  - 5 & - & \ln\left(4  \De  (\sigma_\la, \sigma_\mu, \sigma_\nu)^2 \right).\eeqs
 If
   $\la, \mu, \nu$ are respectively $\sigma_\la \la^{\mathrm{cl}}, \sigma_\mu \la^{\mathrm{cl}}, \sigma_\nu \la^\mathrm{cl}$, then we have equality in (\ref{eq:@}) by Theorem~\ref{thm:5}. This completes the proof of the Theorem.
\end{proof}


Let $h_{0,1,1}^{\mathrm{gue}}$ be the unique continuum hive with the boundary conditions
 $\la = \sigma_\la \lambda^{\mathrm{cl}}$, $\mu = \sigma_\mu \lambda^{\mathrm{cl}}$ and $\nu = \sigma_\nu \lambda^{\mathrm{cl}}$,  for $\sigma_\la = 0, \sigma_\mu = 1, \sigma_\nu = 1$. We define $h_{1,0,1}^{\mathrm{gue}}$ and $h_{1,1,0}^{\mathrm{gue}}$ analogously. 
Given arbitrary $\sigma_\la, \sigma_\mu, \sigma_\nu \in \R$, such that we can find a Euclidean triangle with sides $\sigma_\la, \sigma_\mu, \sigma_\nu$ and area $\De(\sigma_\la, \sigma_\mu, \sigma_\nu) > 0$, define $h_{\sigma_\la, \sigma_\mu, \sigma_\nu}^\mathrm{gue}$ to be the unique continuum hive with spectral boundary conditions $\la = \sigma_\la \lambda^{\mathrm{cl}}$, $\mu = \sigma_\mu \lambda^{\mathrm{cl}}$ and $\nu = \sigma_\nu \lambda^{\mathrm{cl}}$ that is expressible as a (positive) linear combination of $h_{0,1,1}^{\mathrm{gue}}$, $h_{1,0,1}^{\mathrm{gue}}$ and $h_{1,1,0}^{\mathrm{gue}}$.\\

\noindent {\bf Theorem~\ref{thm:9} (extended version).}\\ 
Let $$\varUpsilon(\sigma_\la, \sigma_\mu, \sigma_\nu) := $$ $$\int_T \sigma((-1)(\hess h_{\sigma_\la, \sigma_\mu, \sigma_\nu}^\mathrm{gue})_{ac})\Leb_2(dx) - \inf\limits_{h \in H(\sigma_\la \la^{\mathrm{cl}}, \sigma_\mu \la^{\mathrm{cl}},\la; \sigma_\nu \la^{\mathrm{cl}})}  \int_T \sigma((-1)(\hess h)_{ac})\Leb_2(dx).$$ For $s= (s_0, s_1, s_2) = 2(-\sigma_\la + \sigma_\mu + \sigma_\nu, \sigma_\la - \sigma_\mu + \sigma_\nu, \sigma_\la + \sigma_\mu - \sigma_\nu)\in \R_+^3$, we have
\begin{eqnarray*}  \exp(\frac{5}{4}-\varUpsilon)\frac{(s_0 + s_1 + s_2)s_0s_1s_2}{(s_0 + s_1)(s_1 + s_2)(s_2 + s_0)} \leq  \exp(-\sigma(s))  \leq    \frac{\exp(5)(s_0 + s_1 + s_2)s_0s_1s_2}{36(s_0 + s_1)(s_1+s_2)(s_2 + s_0)}.\end{eqnarray*}

\begin{proof}  
By Green's Theorem and the fact that the largest GUE eigenvalue is asymptotic to $2\sqrt{n}$, we see that the negative of the average Hessian $s$ of $h_{\sigma_\la, \sigma_\mu, \sigma_\nu}^\mathrm{gue}$ satisfies $$\int (-1)(\hess h_{\sigma_\la, \sigma_\mu, \sigma_\nu}^\mathrm{gue})_{ac}\Leb_2(dx) = (s_0, s_1, s_2)$$ $$ = 4(\frac{-\sigma_\la + \sigma_\mu + \sigma_\nu}{2}, \frac{\sigma_\la - \sigma_\mu + \sigma_\nu}{2}, \frac{\sigma_\la + \sigma_\mu - \sigma_\nu}{2}).$$

By the convexity of $\sigma$ and the definition of $h_{\sigma_\la, \sigma_\mu, \sigma_\nu}^\mathrm{gue}$,  we see that 
\beqs \int \sigma((-1)(\hess h_{\sigma_\la, \sigma_\mu, \sigma_\nu}^\mathrm{gue})_{ac})\Leb_2(dx) & \geq &  \sigma( \int (-1)(\hess h_{\sigma_\la, \sigma_\mu, \sigma_\nu}^\mathrm{gue})_{ac}\Leb_2(dx)).\eeqs
Note that by Theorem~\ref{thm:6}, for all $\sigma_\la, \sigma_\mu, \sigma_\nu \in \R$ such that $\De(\sigma_\la, \sigma_\mu, \sigma_\nu) > 0$, we have
\beq \lab{eq:Upsilon}  \int_T \sigma((-1)(\hess h_{\sigma_\la, \sigma_\mu, \sigma_\nu}^\mathrm{gue})_{ac})\Leb_2(dx) - \varUpsilon & =  & \inf\limits_{h \in H(\la, \mu; \nu)}  \int_T \sigma((-1)(\hess h)_{ac})\Leb_2(dx)\nonumber \\ & = & -\frac{5}{4} - \ln\left(\frac{4  \De  (\sigma_\la, \sigma_\mu, \sigma_\nu)^2}{\sigma_\la \sigma_\mu \sigma_\nu}\right).\eeq
From (\ref{eq:Upsilon}) we have 
$\exp( - \sigma(s)) \geq \left(\exp(\frac{5}{4})\frac{(s_0 + s_1 + s_2)s_0s_1s_2}{(s_0 + s_1)(s_1 + s_2)(s_2 + s_0)}\right)\exp(-\varUpsilon)$. 

We now consider a quadratic continuum hive $h^q$ (\ie a quadratic function from $T$ to $\R$ that is also a continuum hive) with  boundary conditions in which the $L^2$ norm of the gradient of the function on three pieces of the boundary are respectively $\sigma_\la, \sigma_\mu$ and $\sigma_\nu$.  Such boundary conditions always have an extension to the interior in the form of a (unique) quadratic continuum hive because as $(\sigma_\la, \sigma_\mu, \sigma_\nu)$ range across this set, all triples satisfying the triangle inequality are attained, but no other triples are attained. Noting that the $L^2$ norm of the identity function $x \mapsto x$ on $[-\frac{1}{2}, \frac{1}{2}]$ equals $\frac{1}{\sqrt{12}}$, we see from (\ref{eq:equality1}) that for the corresponding  (quadratic) boundary values,  $V(\la), V(\mu)$ and $V(\nu)$ are respectively $\sqrt{12} \sigma_\la, \sqrt{12} \sigma_\mu$ and $\sqrt{12} \sigma_\nu$ respectively. Further, (redefining $s$ to the setting of a quadratic hive) the negative of the Hessian $s$ of $h^q$ satisfies $$(-1)(\hess h^q)_{ac} = (s_0, s_1, s_2) = \sqrt{12}(\frac{-\sigma_\la + \sigma_\mu + \sigma_\nu}{2}, \frac{\sigma_\la - \sigma_\mu + \sigma_\nu}{2}, \frac{\sigma_\la + \sigma_\mu - \sigma_\nu}{2}).$$ Therefore,
\beq\lab{eq:take_acc}\frac{(s_0 + s_1 + s_2)s_0s_1s_2}{(s_0 + s_1)(s_1 + s_2)(s_2 + s_0)} = 4\sqrt{12}   \De^2(\sigma_\la, \sigma_\mu, \sigma_\nu).\eeq An upper bound on $-\sigma(s)$ follows by applying Theorem~\ref{thm:6}, since this gives the inequality
\begin{eqnarray*}  \sigma(s)  \geq 3 \ln \sqrt{12} - 5  - \ln\left(\frac{4 \De(\sigma_\la, \sigma_\mu, \sigma_\nu)^2}{\sigma_\la \sigma_\mu \sigma_\nu}\right).\end{eqnarray*}

Taking into account (\ref{eq:take_acc}) this results in 

$$\exp(-\sigma(s))  \leq    \frac{\exp(5)(s_0 + s_1 + s_2)s_0s_1s_2}{36(s_0 + s_1)(s_1+s_2)(s_2 + s_0)}.$$

\end{proof}

Putting aside $\varUpsilon$, the constant in the lower bound on $\exp(-\sigma)$ in Theorem~\ref{thm:9}, namely $\exp(\frac{5}{4})$ is roughly 
$3.49034296$, while the constant $\frac{\exp(5)}{36}$ in the upper bound on $\exp(-\sigma)$ in Theorem~\ref{thm:9} is roughly $4.12258775$.
 As Figure~\ref{fig:GUE_3,4,5} shows, a non-rigorous heuristic leads us to believe that at least for values of $s_0, s_1, s_2$, that are nearly equal, $\varUpsilon$ is close to $0$.

 \section{Simulations of random hives and  random lozenge tilings via the octahedron recurrence}\label{recurrence-sec}

   \begin{figure}
\begin{center}
\includegraphics[scale=0.40]{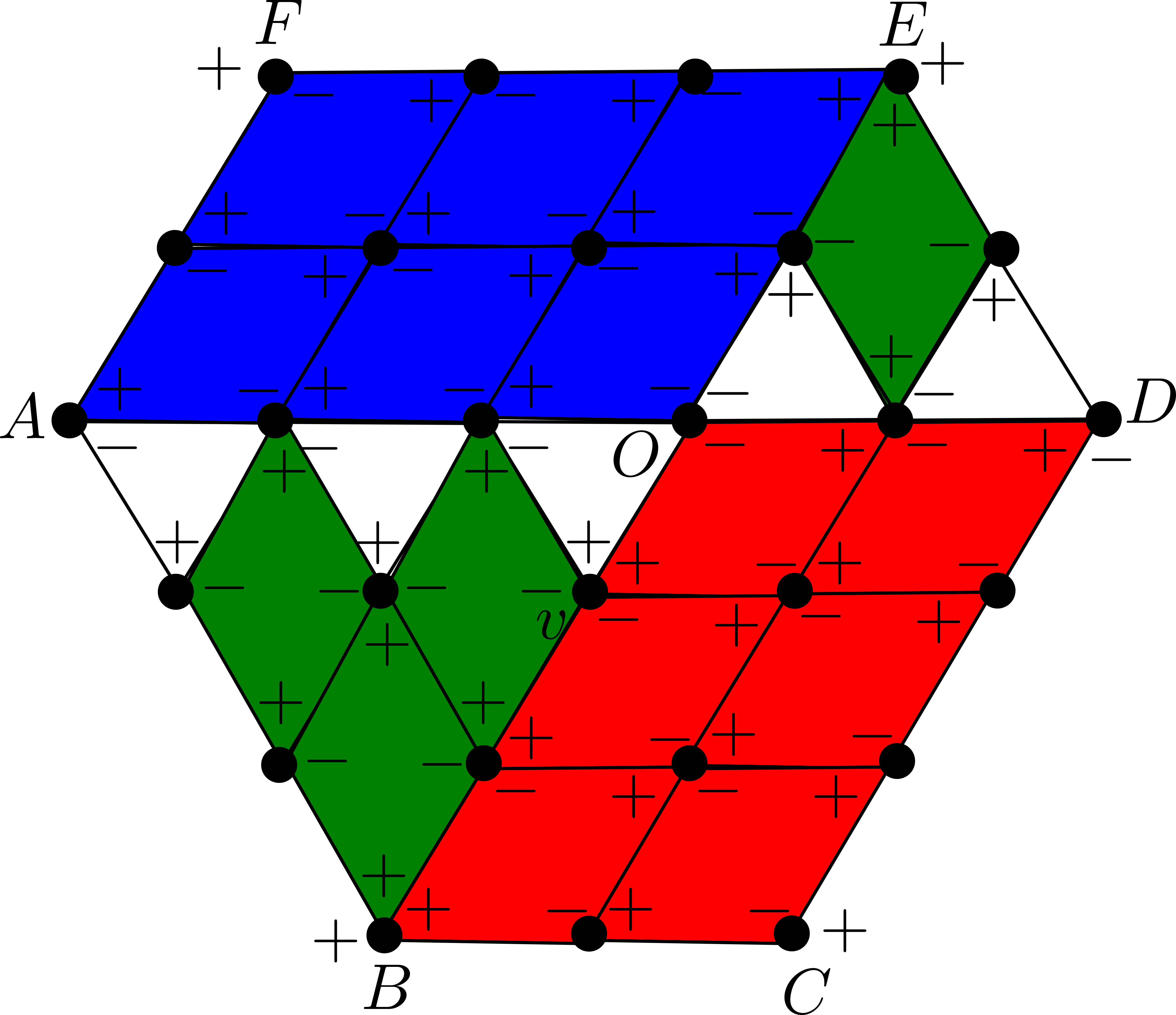}
\caption{The standard lozenge tiling of a hexagon $ABCDEF$ centered at $v$.  The total weight of this tiling is $\tilde k(E) + \tilde k(B) - \tilde k(O)$, where $O$ is the intersection of the diagonal $BE$ with the equator $AD$.}\label{fig:standard}
\end{center}
\end{figure}

\begin{figure}
\begin{center}
\includegraphics[scale=0.40]{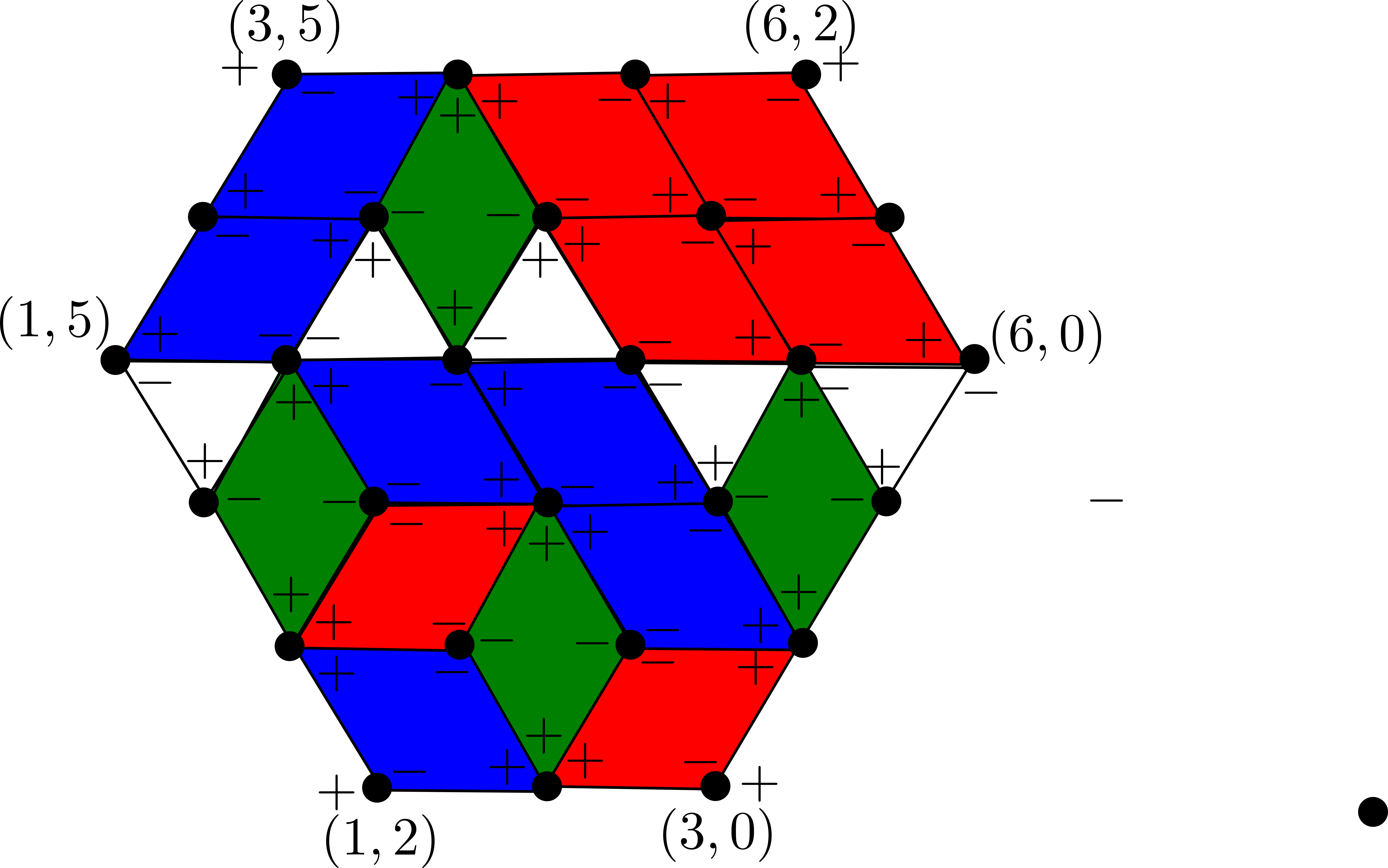}
\caption{A typical lozenge tiling of $\hexagon_{(3,2)}^6$, $n=6$. }\label{fig:typical}
\end{center}
\end{figure}

\begin{figure}[!ht]
  \centering

  \begin{subfigure}{0.15\linewidth}
    \includegraphics[width=\linewidth]{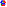}
    \caption{$\hexagon_{(2,2)}^{4}$, $n=4$}
    \label{fig:loz:4}
  \end{subfigure}\hspace{1em}
  \begin{subfigure}{0.2\linewidth}
    \includegraphics[width=\linewidth]{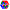}
    \caption{$\hexagon_{(3,3)}^{6}$, $n=6$}
    \label{fig:loz:6}
  \end{subfigure}\hspace{1em}
  \begin{subfigure}{0.3\linewidth}
    \includegraphics[width=\linewidth]{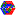}
    \caption{$\hexagon_{(5,5)}^{10}$, $n=10$}
    \label{fig:loz:10}
  \end{subfigure}

  \vspace{1em}

  \begin{subfigure}{0.5\linewidth}
    \includegraphics[width=\linewidth]{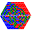}
    \caption{$\hexagon_{(10,10)}^{20}$, $n=20$}
    \label{fig:loz:20}
  \end{subfigure}

  \caption{Random lozenge tilings derived via the octahedron recurrence from pairs of hives $(k,k')$ sampled from GUE minor processes.}
  \label{fig:lozenge_A}
\end{figure}

\begin{figure}[!ht]
  \centering



  \begin{subfigure}{0.5\linewidth}
    \includegraphics[width=\linewidth]{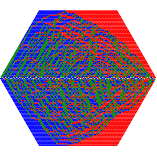}
    \caption{$\hexagon_{(50,50)}^{100}$, $n=100$}
    \label{fig:loz:100:1}
  \end{subfigure}

  \vspace{1em}

  \begin{subfigure}{0.5\linewidth}
    \includegraphics[width=\linewidth]{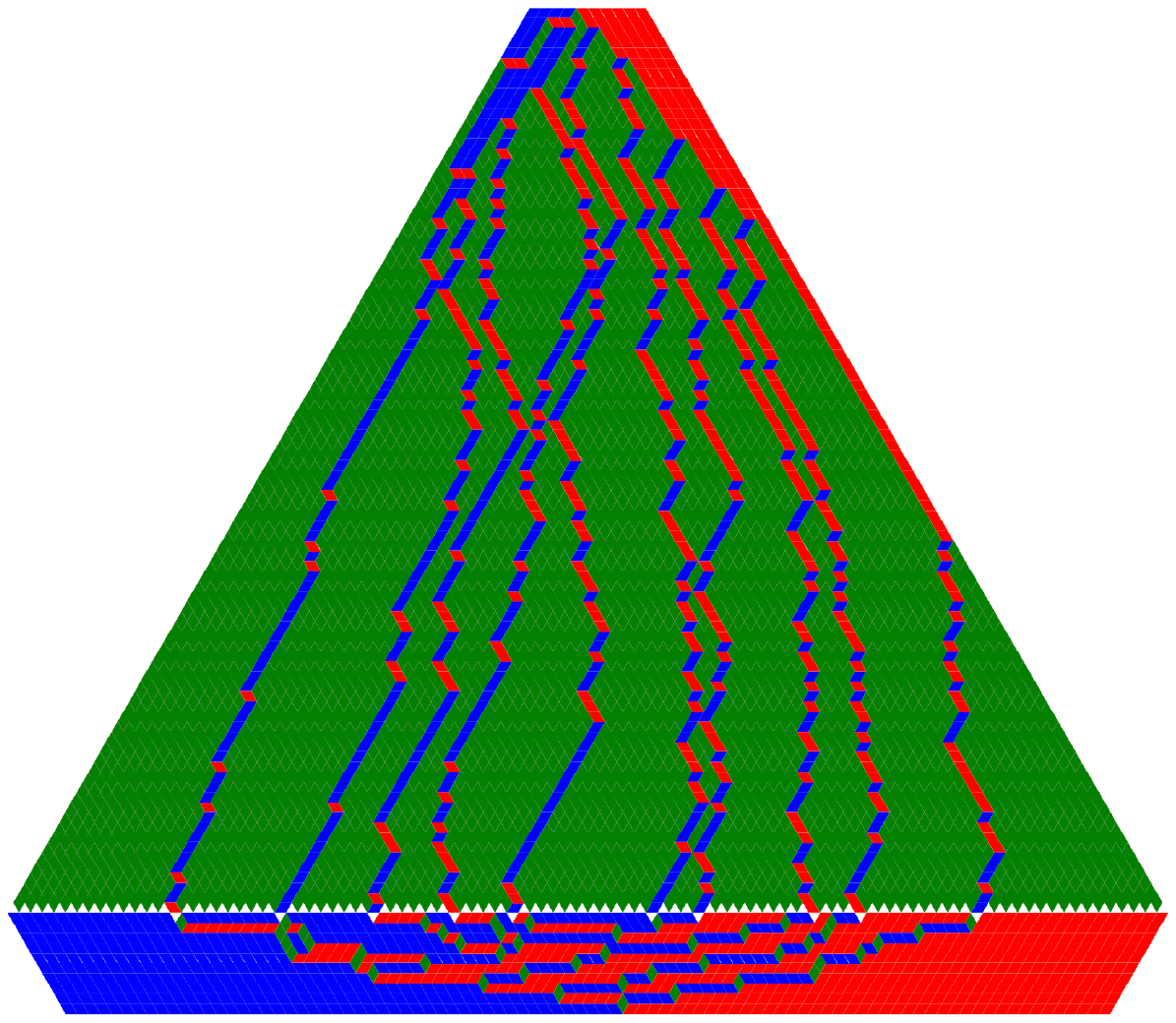}
    \caption{$\hexagon_{(90,90)}^{100}$, $n=100$}
    \label{fig:loz:100:skew}
  \end{subfigure}

  \caption{Random lozenge tilings as in Fig.~\ref{fig:lozenge_A}. 
  }
  \label{fig:lozenge_B}
\end{figure}

We refer the reader to  \cite[\S 3]{NarSheffTao} for further background, but reproduce some of the material from there that we will need to explain the simulations Figure~\ref{fig:lozenge_A} and Figure~\ref{fig:lozenge_B}.

For the purposes of this section only, $T := \{(i,j) : 0 \leq i \leq j \leq n\}$ and $T' := \{(i,j) : 0 \leq j \leq i \leq n\}$ and  $U := \{(i,j) : i,j \leq n \leq i+ j\}$ and $U' := \{(i, j):i, j \leq i + j \leq n\}.$

\begin{definition}[Lozenges and border triangles]  A \emph{lozenge} is a quadruple $ABCD$ in $U$ or $U'$ that is one of following three forms for some $i,j \in \Z$: (For the coordinates, see Figure~\ref{fig:typical}.)
\begin{itemize}
\item[(i)] $(A,B,C,D) = ((i,j), (i+1,j-1), (i+2,j-1), (i+1,j))$
\item[(ii)] $(A,B,C,D) = ((i,j), (i,j+1), (i-1,j+2), (i-1,j+1))$
\item[(iii)] $(A,B,C,D) = ((i,j), (i+1,j), (i+1,j+1), (i,j+1))$.
\end{itemize}
Lozenges of type (i) will be called \emph{blue} if they lie in the upper triangle $U$ and \emph{red} if they lie in the lower triangle $U'$; lozenges of type (ii) will be called \emph{red} if they lie in $U$ and \emph{blue} if they lie in $U'$; and lozenges of type (iii) that lie either in $U$ or in $U'$ will be called \emph{green}; see Figure \ref{fig:typical}.  A quadruple of the form (iii) that crosses the diagonal separating $U$ and $U'$ is \emph{not} considered to be a lozenge, but instead splits into two border triangles as defined below.  (The colors of lozenges will not be needed immediately, but will play a useful role later in this section.)

A \emph{border edge} is an edge $AC$ of the form $(A,C) = ((i,n-i), (i+1,n-i-1))$ for some $0 \leq i < n$; the border edges thus separate $U$ and $U'$.  Each border edge $(A,C) = ((i,n-i), (i+1,n-i-1))$ is bordered by two \emph{border triangles} $ABC$, defined as follows:
\begin{itemize}
\item (Upward triangle) $(A,B,C) = ((i,n-i), (i+1,n-i), (i+1,n-i-1))$.
\item (Downward triangle) $(A,B,C) = ((i,n-i), (i,n-i-1), (i+1,n-i-1))$.
\end{itemize}
Note that upward triangles lie (barely) in $U$, while downward triangles lie (barely) in $U'$.

Given a lozenge $\edge = ABCD$ and a function $\tilde k \colon \{0,\dots,n\}^2 \to \R$ defined as before, we define the \emph{weight} $\weight(\edge) = \weight(\edge, \tilde k)$ to be the quantity
$$\weight(\edge) \coloneqq \frac{1}{3} (\tilde k(A) + \tilde k(C) - \tilde k(B) - \tilde k(D)).$$
Similarly, given a border triangle $\Delta = ABC$, the weight $\weight(\Delta) = \weight(\Delta,  \tilde k)$ is defined as
$$\weight(\Delta) \coloneqq \frac{1}{3} (\tilde k(B) -  \tilde k(A)).$$ \end{definition}

\begin{definition}[Octahedron recurrence]\label{octa-def}  If $v = (i,j)$ lies in the interior of $\{0,\dots,n\}^2 = T \cup T'$, then the \emph{excavation hexagon}
$\hexagon_v^n = ABCDEF$ in $\{0,\dots,n\}^2 = U \cup U'$ centered at $v$ is defined as follows:
\begin{itemize}
\item If $v \in T$ (i.e., $i \leq j$), then
$$ (A,B,C,D,E,F) = ((0,n), (0,j), (i,j-i), (n+i-j,j-i), (n+i-j,j), (i,n)).$$
\item If $v \in T'$ (i.e., $i \geq j$), then
$$ (A,B,C,D,E,F) = ((i-j, n+j-i), (i-j,j), (i,0), (n,0), (n,j), (i,n+j-i)).$$
\end{itemize}
Note that these two definitions agree when $v \in T \cap T'$ (i.e., when $i=j$). The original point $v = (i,j)$ is then the intersection of the diagonals $BE$ and $CF$.  The line $AD$ is called the \emph{equator}; it lies on the border between $U$ and $U'$.  The \emph{weight} $\weight(\hexagon_v^n) = \weight(\hexagon_v^n,\tilde k)$ of this hexagon is defined as
\begin{equation}\label{hexagon-weight}
 \weight(\hexagon_v^n) \coloneqq \frac{1}{3} (\tilde k(B) + \tilde k(C) - \tilde k(D) + \tilde k(E) + \tilde k(F)).
\end{equation}

A \emph{lozenge tiling} $\Xi$ of the excavation hexagon $\hexagon_v^n$ is a partition of the (solid) hexagon into (solid) lozenges and (solid) border triangles, such that each border edge on the equator is adjacent to exactly one border triangle in the tiling; see Figure \ref{fig:typical}.  An example of a lozenge tiling is the \emph{standard lozenge tiling} $\Xi_0$, in which the trapezoid $ABEF$ is tiled by blue lozenges in $U$ and by green lozenges and downward border triangles in $U'$, while the opposite trapezoid $BCDE$ is tiled by green lozenges and upward border triangles in $U$ and by red lozenges in $U'$; see Figure \ref{fig:standard}.

The \emph{weight} $w_\Xi = w_\Xi(\tilde k)$ of such a tiling is defined to be the sum of the weights of all the lozenges $\edge$ and triangles $\Delta$ in the tiling, as well as the weight of the entire hexagon $\hexagon_v^n$:
\begin{equation}\label{weight-form}
w_\Xi \coloneqq \sum_{\edge \in \Xi} \weight(\edge) + \sum_{\Delta \in \Xi} \weight(\Delta) + \weight(\hexagon_v^n).
\end{equation}
Note that the $w_\Xi$ depend linearly on $\tilde k$, and hence on $k, k'$.  We then define
\beq\lab{eq:hv} \tilde h(v) \coloneqq \max_{\Xi\ \mathrm{tiles}\ \hexagon_v^n} w_\Xi.\eeq
\end{definition}

As remarked on  \cite[page 15]{NarSheffTao},  using Speyer's Theorem \cite{Speyer} and the volume preserving nature of the octahedron recurrence, it can be seen that when $k$ and $k'$ are hives arising via Proposition~2(iv) of \cite{NarSheffTao} from random Gelfand-Tsetlin patterns sampled independently from two scaled GUE minor processes, the distribution of $\tilde h(v)$ is exactly the same as the distribution of the value of the random augmented hive at $v$.

It was shown in \cite{NarSheffTao} that we can also write the weight form $w_\Xi$ in \eqref{weight-form} in a red lozenge-avoiding form as
\begin{equation}\label{alt-weight}
w_\Xi \coloneqq 2 \sum_{\edge \in \Xi\text{, blue}} \weight(\edge) + \sum_{\edge \in \Xi\text{, green}} \weight(\edge) 
+ \sum_{\Delta \in \Xi} \weight(\Delta) + \weight'(\hexagon_v^n)
\end{equation}
where the modified weight $\weight'(\hexagon_v)$ of the hexagon $\hexagon_v$ is defined by the formula
$$
 \weight'(\hexagon_v^n) \coloneqq \frac{1}{3} (-\tilde k(A) + 2\tilde k(B) + 2\tilde k(F)).
$$
 \begin{remark}[Role of large gaps]\lab{rem:2}  As just mentioned, the red lozenges can be completely avoided, and since it is only the red lozenges whose weights depend on the large gaps, this shows that as long as it is Speyer's Theorem that is used, and not the iterative octahedron recurrence, then for $v \in T$, the value $\tilde h(v)$ is independent of the specific values of the gaps. This in turn implies that even if one uses the iterative octahedron recurrence, as long as the gaps are large enough, for $v \in T$, $\tilde h(v)$ is independent of them. The lozenge tilings that are shown in Figure~\ref{fig:lozenge_A} and \ref{fig:lozenge_B} were obtained by computing maximum weight matchings for the weighted graph as in (\ref{eq:hv}) by optimizing over a certain matching polytope.
 \end{remark}

\begin{figure}
  \begin{center}
  \includegraphics[scale=0.40]{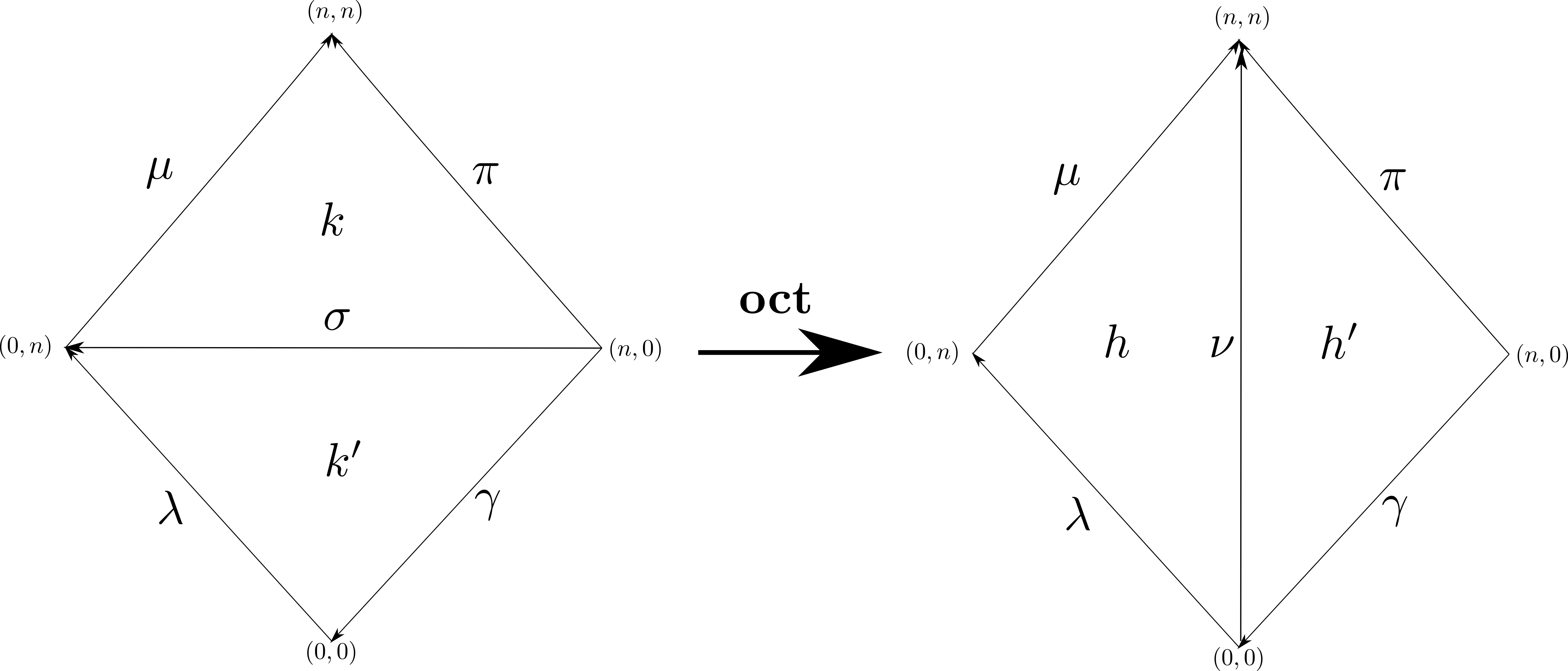}
  \caption{A schematic depiction from \cite{NarSheffTao} of the octahedron recurrence that transforms one pair $(k,k') $ of hives into another $(h,h')$.  The hives $h,h',k,k'$ have been shifted to lie on triangles $T, T', U, U'$ respectively. }\label{fig:octahedron}
  \end{center}
  \end{figure}

\begin{figure}[!h]
    \centering
    \includegraphics[width=1\linewidth]{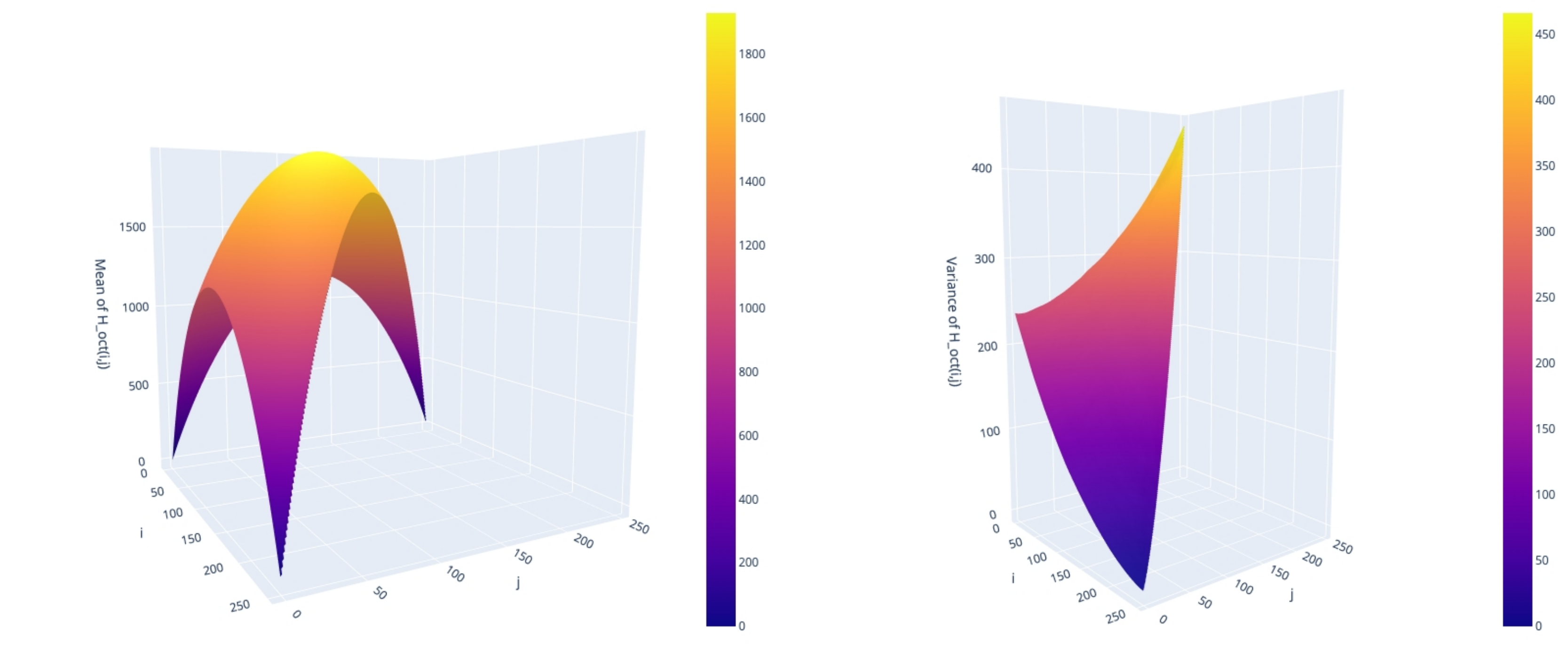}
    \caption{The pointwise mean (image on the left) and pointwise variance (image on the right) of the bottom panel hive after octahedron recurrence from minor processes obtained from two GUE matrices with standard deviation $1$.}
    \label{fig:GUE_Octahedron_Recurrence_1_1}
\end{figure}

\begin{figure}[!h]
    \centering
    \includegraphics[width=1\linewidth]{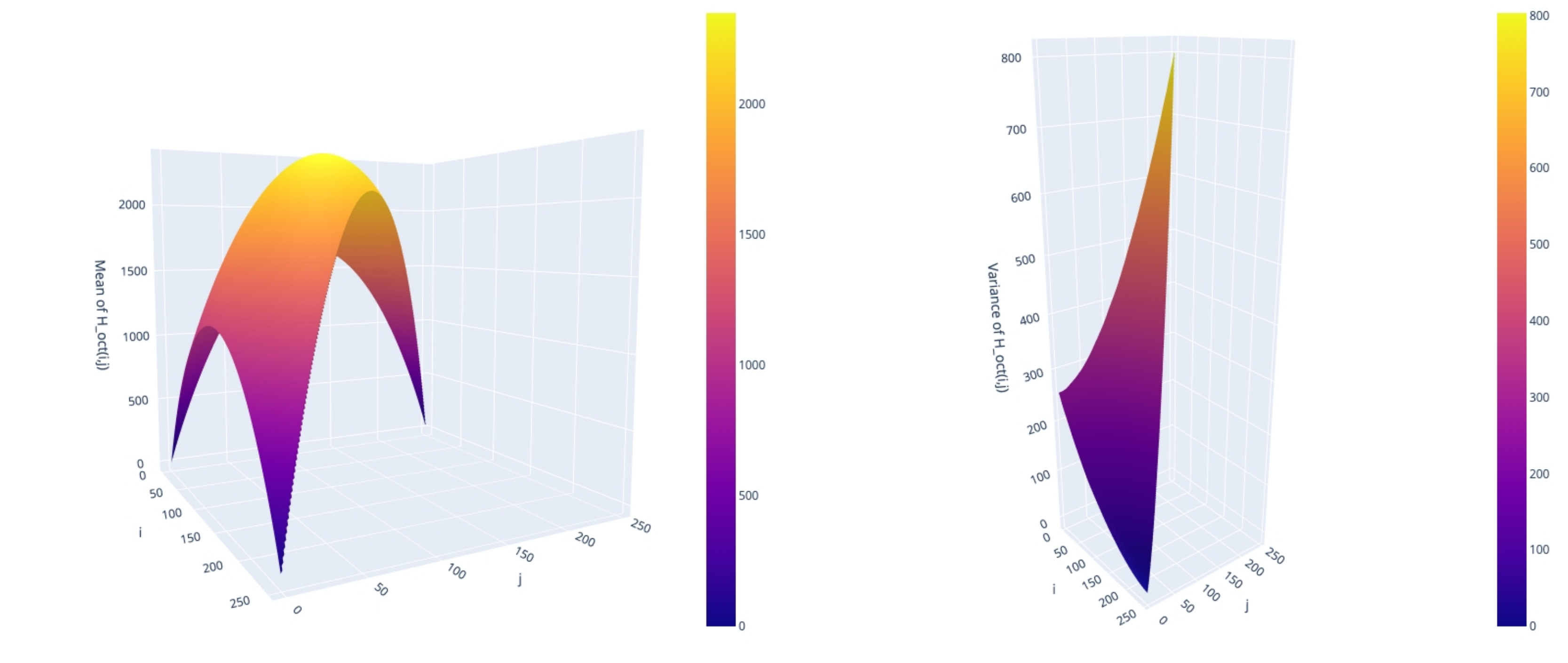}
    \caption{The pointwise mean (image on the left) and pointwise variance (image on the right) of the bottom panel hive after octahedron recurrence from minor processes obtained from two GUE matrices with standard deviation $1$ and $\sqrt{2}$ respectively.}
    \label{fig:GUE_Octahedron_Recurrence_1_sqrt2}
\end{figure}

\begin{figure}[!h]
    \centering
    \includegraphics[width=1\linewidth]{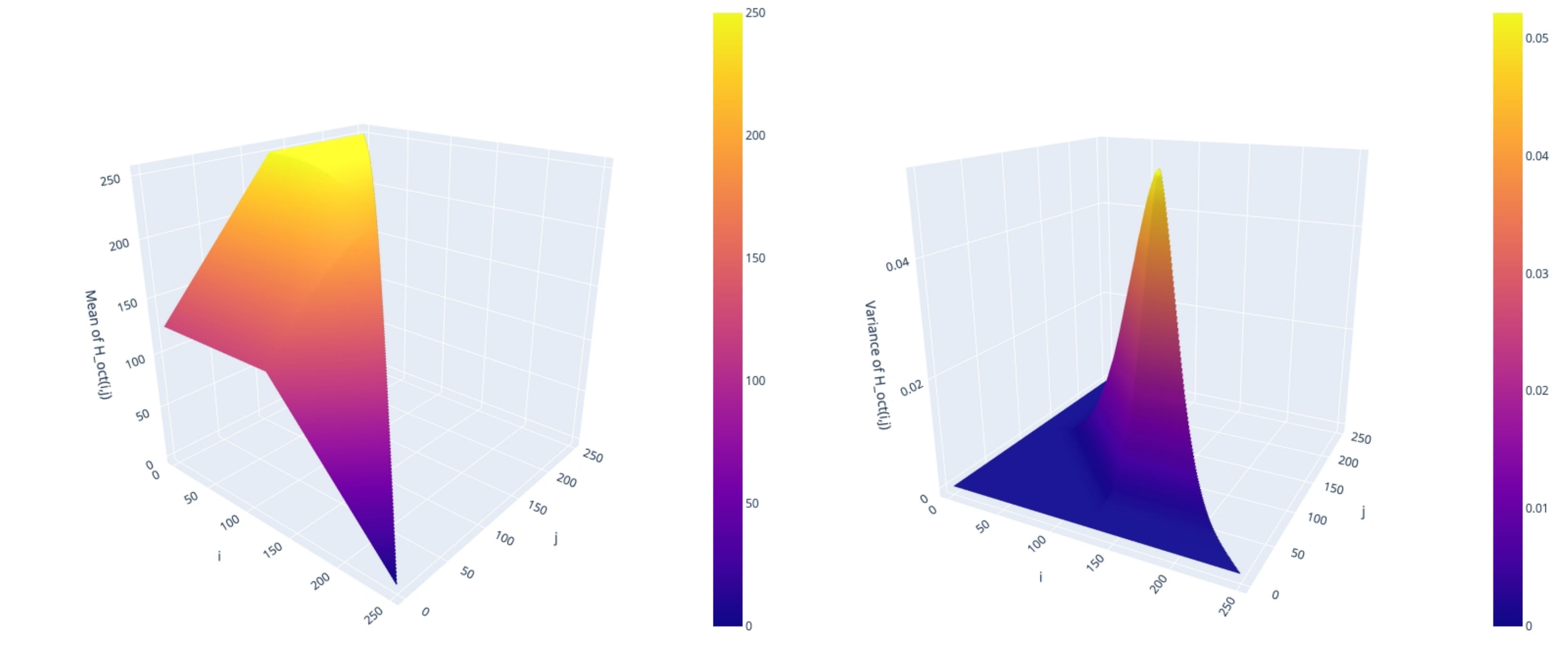}
    \caption{The pointwise mean (image on the left) and pointwise variance (image on the right) of the bottom panel hive after octahedron recurrence from minor processes obtained from two random projections of rank $\frac{n}{2}$ from unitarily invariant distributions.}
    \label{fig:RP_Octahedron_Recurrence}
\end{figure}

\begin{observation}
In the lower trapezoid (contained in $k'$ in Figure~\ref{fig:octahedron}),
for a green lozenge $\edge$ that has labels $\la_{i, j}$ and $\la_{i, j+1}$ on two opposite sides, $$\wt(\edge) :=  \frac{1}{3}\left(\la_{i, j} - \la_{i, j+1}\right),$$ and a blue lozenge  $\edge$ that has labels $\la_{i, j}$ and $\la_{i-1, j-1}$ on two opposite sides, $$\wt(\edge) :=  \frac{1}{3}\left(\la_{i, j} - \la_{i-1, j-1}\right).$$
In the upper trapezoid (contained in $k$ in Figure~\ref{fig:octahedron}),
for a blue lozenge $\edge$ that has labels $\mu_{i, j}$ and $\mu_{i, j+1}$ on two opposite sides, $$\wt(\edge) :=  \frac{1}{3}\left(\mu_{i, j} - \mu_{i, j+1}\right),$$ and a green lozenge  $\edge$ that has labels $\mu_{i, j}$ and $\mu_{i-1, j-1}$ on two opposite sides, $$\wt(\edge) :=  \frac{1}{3}\left(\mu_{i, j} - \mu_{i-1, j-1}\right).$$
\end{observation}

\section{Estimating the surface tension function of hives numerically}
\label{sec:EstimatingSurfaceTension}

In this section, we estimate the surface tension function of hives numerically.

Recall from Theorem~\ref{thm:6} that the surface tension function $\sigma$ satisfies the following equation when $\sigma_\la, \sigma_\mu, \sigma_\nu$ form a triangle with $\De  (\sigma_\la, \sigma_\mu, \sigma_\nu) > 0$:
\begin{eqnarray} \lab{eq:sig}\inf\limits_{h \in H(\la, \mu; \nu)}  \int_T \sigma((-1)(\hess h)_{ac})\Leb_2(dx) = -\frac{5}{4} - \ln\left(\frac{4  \De  (\sigma_\la, \sigma_\mu, \sigma_\nu)^2}{\sigma_\la \sigma_\mu \sigma_\nu}\right). \end{eqnarray}

Thus, if we knew the shape of a surface tension minimizer for all GUE boundary conditions, finding $\sigma$ potentially reduces to solving an infinite family of linear equations. 

Though it would be natural to expect that any subsequential limit of the center of masses of the hive polytopes as $n \ra \infty$ would be a minimizer,  for reasons of computational efficiency, we instead use the following heuristic approach. For this reason, the results of this section should be taken as merely suggestive.

\begin{figure}[!h]
    \centering
    \includegraphics[width=0.5\linewidth]{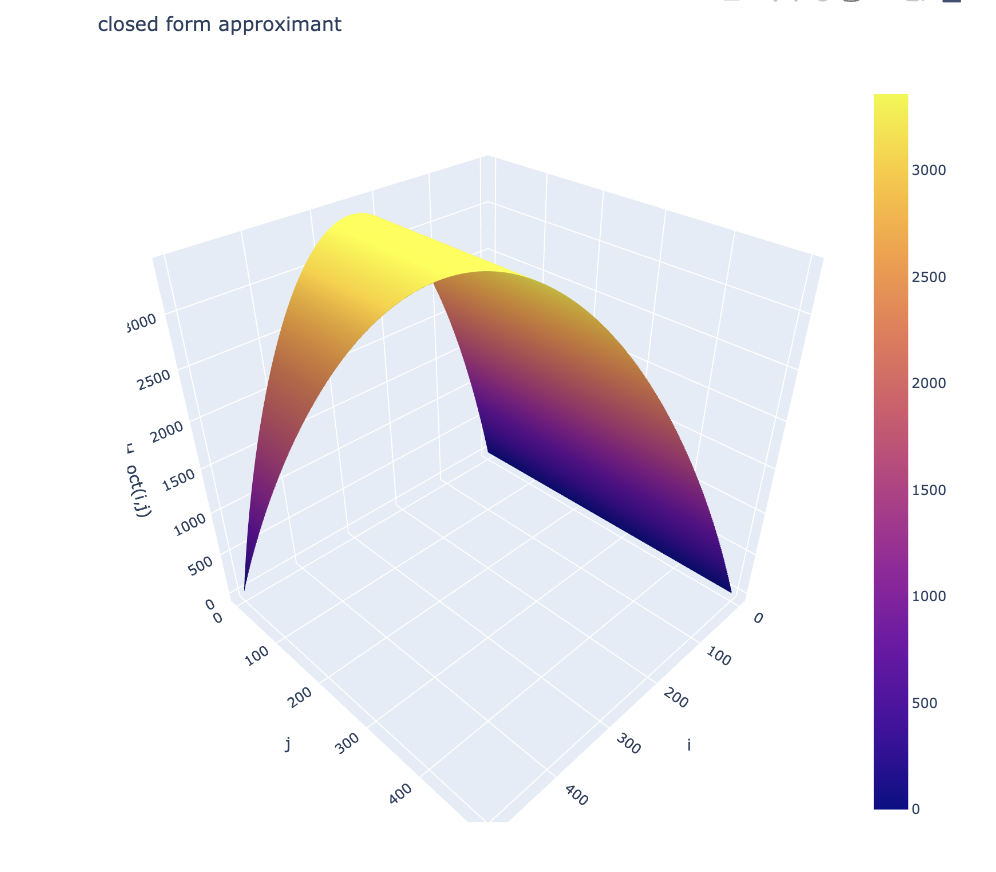}
    \caption{GUE boundary conditions with $\sigma_\la, \sigma_\mu, \sigma_\nu$ respectively equal to $0, 1, 1$ for $n=500$.
    Positive linear combinations of three such hives are used to produce GUE hives with arbitrarily scaled GUE boundary conditions.}
    \label{fig:GUE_0,1,1}
\end{figure}
\newpage
\begin{figure}[!h]
    \centering
    \includegraphics[width=0.5\linewidth]{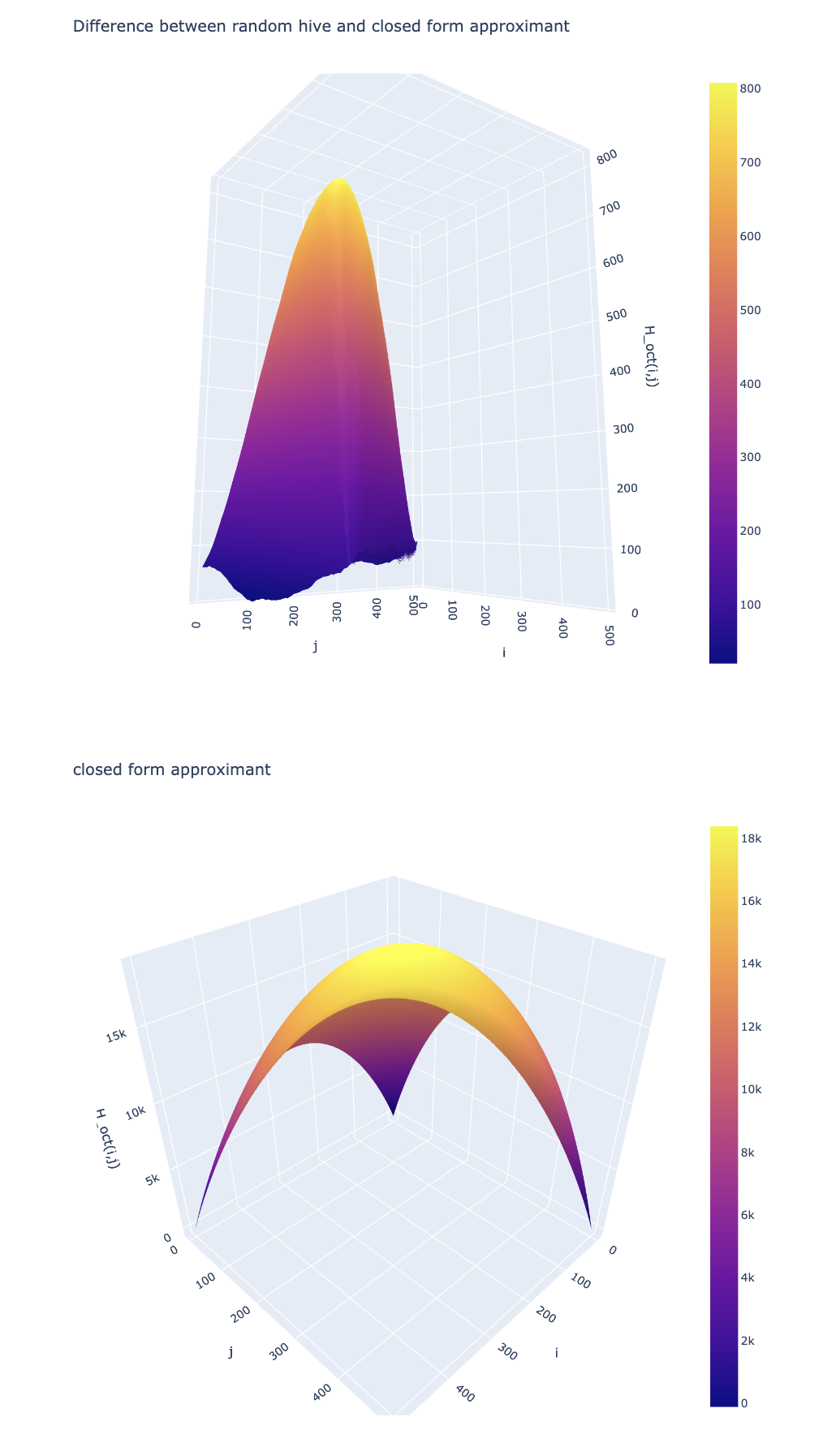}
    \caption{A random hive in the lower figure with GUE boundary conditions with $\sigma_\la, \sigma_\mu, \sigma_\nu$ respectively equal to $3, 4, 5$ for $n=500$. In the upper figure, we subtract from the random hive an interpolant that is a suitable positive linear combination of the three rotated copies of the hive in  Figure~\ref{fig:GUE_0,1,1}.
    }
    \label{fig:GUE_3,4,5}
\end{figure}

\begin{figure}[!h]
    \centering
    \includegraphics[width=0.5\linewidth]{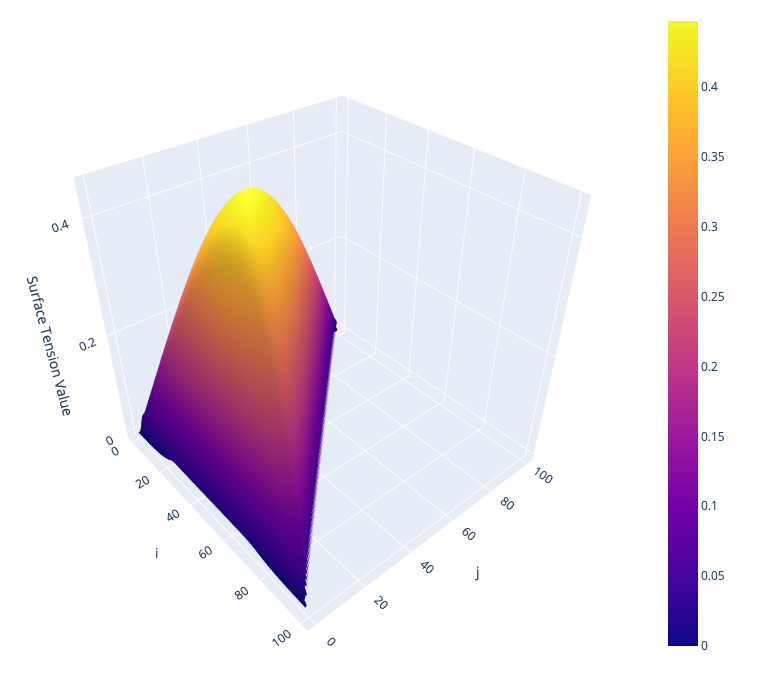}
    \caption{A plot of a numerical evaluation of $\exp(-\sigma^{\mathtt{num}}(s))$ where $s$  lies on the triangle $\{s \in \R_+^3: s_0 + s_1 + s_2 = 1\}$. }
    \label{fig:exp-sig}
\end{figure}

\begin{figure}[!h]
    \centering
    \includegraphics[width=0.5\linewidth]{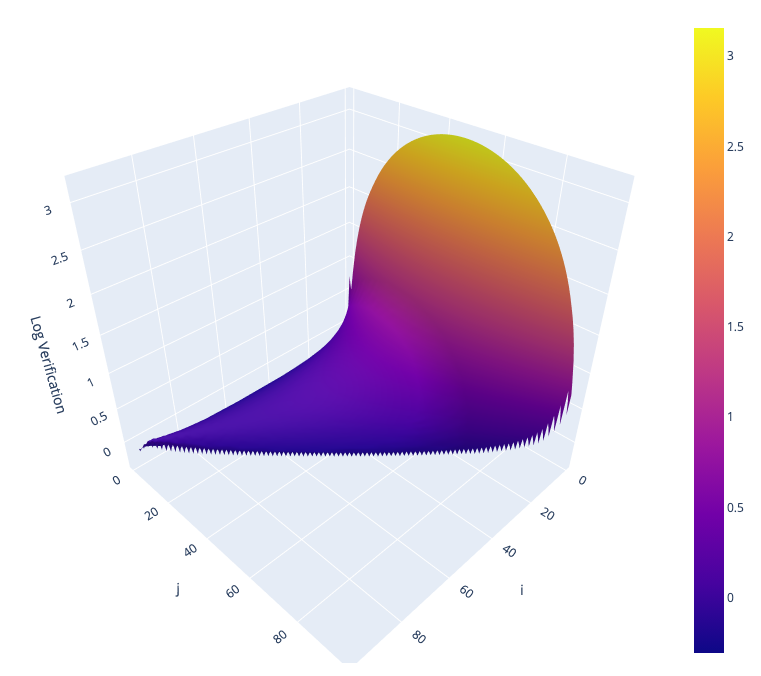}
    \caption{A plot of the numerical evaluation of $\sigma^{\mathtt{num}}(s) - f(s)$ where $s$  lies on the triangle $\{s \in \R_+^3: s_0 + s_1 + s_2 = 1\}$ and $f(s) = - \ln(\frac{e}{\pi} (s_0 + s_1) \sin(\frac{\pi s_0}{s_0 + s_1}))$. $f(s)$ is the surface tension when $s_2 \ra \infty.$ 
    }
    \label{fig:sig-f}
\end{figure}\bibliographystyle{alpha}

\begin{figure}[!h]
    \centering
    \includegraphics[width=0.5\linewidth]{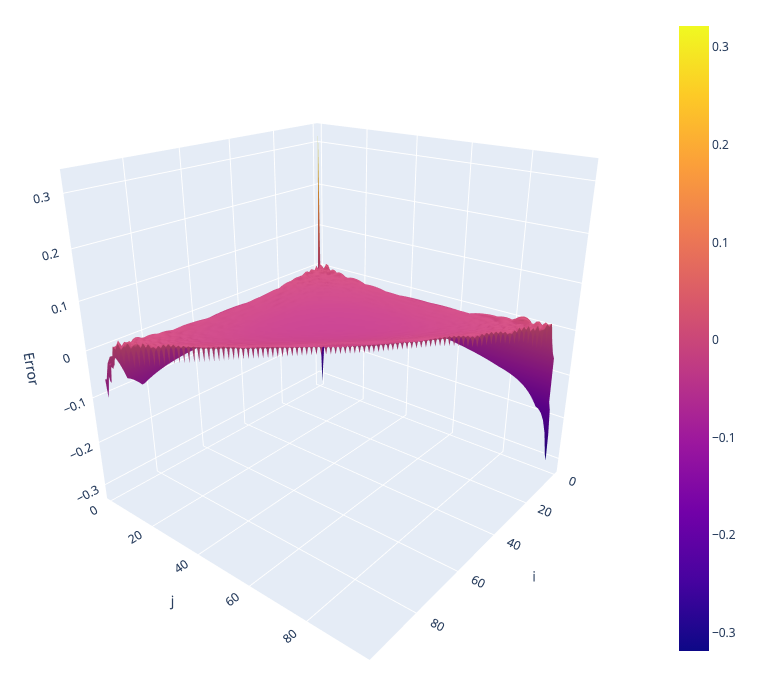}
    \caption{A plot of $\sigma^{\mathtt{num}}(s) - g(s)$ where $g(s) = -\ln\left(\frac{4e(s_0 + s_1 + s_2) s_0 s_1 s_2}{\pi(s_0 + s_1)(s_1 + s_2) (s_2 + s_0)}\right)$ and  $s$  lies on the triangle $\{s \in \R_+^3: s_0 + s_1 + s_2 = 1\}$. 
    .}
    \label{fig:sig-estimate}
\end{figure}


We consider GUE boundary conditions with $\sigma_\la, \sigma_\mu, \sigma_\nu$ respectively equal to $0, 1, 1$, (which has a unique hive extension) and consider positive linear combinations of these unique extensions to give a hive with arbitrary scaled GUE boundary conditions. Note that when we subtract from the random hive we generated with GUE boundary conditions with $\sigma_\la, \sigma_\mu, \sigma_\nu$ respectively equal to $3, 4, 5$, such an interpolant that is a suitable positive linear combination of the three rotated copies of the hive in Figure~\ref{fig:GUE_0,1,1},  we see when $n = 500$ in Figure~\ref{fig:GUE_3,4,5}, that the $L^\infty$ norm of the error is about $800$, which is roughly $5$ percent of the $L^\infty$ norm of the interpolant, which is about $16000$. As mentioned in the end of Section~\ref{sec:bounds_surf}, this empirical fact, suggests in a non-rigorous fashion that $\varUpsilon(3, 4, 5)$ is small in magnitude.

For $n = 100$, we choose $10^4$ random hives (see Section~\ref{recurrence-sec} for the sampling methodology) with GUE boundary conditions with $\sigma_\la, \sigma_\mu, \sigma_\nu$ respectively, and minimize the residual quadratic loss to solve the overdetermined system of 
equations arising from (\ref{eq:sig}) for a discretization of $\sigma$, while enforcing the constraint that $\sigma$ is convex. As a result, we obtained the plot in figure~\ref{fig:exp-sig}.
Suppose that  $f(s)$ is the limiting surface tension of the hive model as $s_2 \ra \infty$. 
 Because the polytope associated with Gelfand-Tsetlin patterns is a limit of the hive model as $s_2 \ra \infty$ (as explained in Proposition~2(iv) of \cite{NarSheffTao}) it follows from \cite[Corollary 1.5]{Johnston} that the surface tension of a Gelfand-Tsetlin pattern is 
$$f(s) = - \ln(\frac{e}{\pi} (s_0 + s_1) \sin(\frac{\pi s_0}{s_0 + s_1})).$$ 
The numerical results in Figures~\ref{fig:exp-sig}, appear to satisfy (see Figure~\ref{fig:sig-f}) \begin{eqnarray}\lab{eq:sig-f} \sigma^{\mathtt{num}}(s) - f(s) \ra 0 \end{eqnarray} as $s_2 \ra \infty$.  This is  in accordance with the results of Shlyakhtenko and Tao \cite{SchTao} and Johnston \cite{Johnston} which predicts the limit to be $0$. 

By \cite[Corollary 2]{NarSheff}, we know that $\exp(-\sigma(s))$ is a concave function of $s$ for $s \in \R_+^3$ and the numerical results in Figure~\ref{fig:exp-sig} are consistent with this.

We find that $\sigma^{\mathtt{num}}(s)$ is close to  $$g(s):= -\ln\left(\frac{4e(s_0 + s_1 + s_2) s_0 s_1 s_2}{\pi(s_0 + s_1)(s_1 + s_2) (s_2 + s_0)}\right).$$ The form of $g(s)$ was taken to be a scalar multiple of the upper and lower bounds in Theorem~\ref{thm:9}. The limit of $g(s) - f(s)$ as $s_2 \ra \infty$ does not equal $0$ for all $s_0$ and $s_1$ (though $g(s)$ has been designed to satisfy this along the line segment where $s_0 = s_1$,) thus $\sigma \neq g$ in general. A numerical comparison between $\exp(-\sigma^{\mathtt{num}}(s))$ and $\exp(-g(s))$ is depicted in Table \ref{tab:exp_sigma_values} at selected values of $s$.

\begin{table}[!h]
\centering
\begin{tabular}{|c|c|c|c|c|}
\hline
\textbf{$s_0$} & \textbf{$s_1$} & \textbf{$s_2$} & \textbf{$\exp(-\sigma^{\mathtt{num}}(s))$} & \textbf{$\exp(-g(s))$} \\ \hline
0.1 & 0.1 & 0.8 & 0.172 & 0.171 \\ \hline
0.1 & 0.2 & 0.7 & 0.227 & 0.224 \\ \hline
0.1 & 0.3 & 0.6 & 0.251 & 0.247 \\ \hline
0.1 & 0.4 & 0.5 & 0.261 & 0.256 \\ \hline
0.2 & 0.2 & 0.6 & 0.332 & 0.324 \\ \hline
0.2 & 0.3 & 0.5 & 0.381 & 0.371 \\ \hline
0.2 & 0.4 & 0.4 & 0.396 & 0.385 \\ \hline
0.3 & 0.3 & 0.4 & 0.438 & 0.424 \\ \hline
\end{tabular}
\caption{Comparison between $\exp(-\sigma^{\mathtt{num}}(s))$ and $\exp(-g(s))$ at selected values of $s=(s_0,s_1,s_2)$.}
\label{tab:exp_sigma_values}
\end{table}

\subsection{Remarks}\lab{sec:remarks}
Let $$ (s_0, s_1, s_2) = 4(\frac{-\sigma_\la + \sigma_\mu + \sigma_\nu}{2}, \frac{\sigma_\la - \sigma_\mu + \sigma_\nu}{2}, \frac{\sigma_\la + \sigma_\mu - \sigma_\nu}{2}).$$ 
The difference between the quantity in the upper bound in Theorem~\ref{thm:9}  and $g(s)$ is $$ - \ln \frac{4e}{\pi} + \frac{5}{4} - \varUpsilon \approx 0.00843552 - \varUpsilon.$$
In general, the sign of $ - \ln \frac{4e}{\pi} + \ln (\exp(5/4)) - \varUpsilon(\sigma_\la, \sigma_\mu, \sigma_\nu)$ is unknown to us. However,  the fact that $g$ (by design) satisfies $$\lim\limits_{s_2 \ra \infty} g(1, 1, s_2) - f(1, 1) = 0,$$ implies that $$\lim\limits_{s_2 \ra \infty} \varUpsilon (\sigma_\la, \sigma_\mu, \sigma_\nu) > 0,$$ where the limit is taken along the line on which $s_0 = s_1$.  This implies that the function $h_{\sigma_\la, \sigma_\mu, \sigma_\nu}^{\mathrm{gue}}$  used in the proof of Theorem~\ref{thm:9} and in our numerical experiments is not in general, a minimizer of the surface tension. 


\section{Acknowledgements}
We are grateful to Terence Tao for encouragement and the suggestion to consider the plot in Figure~\ref{fig:sig-f}. We thank the anonymous referee for the careful reading and helpful suggestions. 
We also acknowledge the use of large language models as an aid in checking exposition and editing parts of the manuscript; all mathematical arguments and any remaining errors are the responsibility of the authors.
 HN acknowledges support from the Department of Atomic Energy, Government of India, under project no. RTI4001, a Swarna Jayanthi fellowship and the Infosys-Chandrasekharan virtual center for Random Geometry.



\appendix

\section{A result of Landon and Sosoe from \cite{Sosoe}}\lab{ap:Sosoe}

We now describe a result of Landon and Sosoe, namely Theorem 1.4 of \cite{Sosoe}.

 Let $\xi_o$ and $\xi_d$ be two real centered random variables with bounded moments of all orders, with the variance of $\xi_o$  being $1$. A real symmetric Wigner matrix is an $N \times N$ self adjoint matrix so that the entries $(H_{ij})_{i\leq j}$ are independent, and
$$\sqrt{N}{H_{ii}} \sim \xi_d,$$ $$ \sqrt{N}H_{ij} \sim \xi_o, i \neq j.$$
Let \( H \) be a real symmetric Wigner matrix, and let \( a_k \) and \( s_k \) be the cumulants of the matrix entries. Denote by \( \xi_o \) and \( \xi_d \) independent standard Gaussian random variables and set
\begin{equation}
\begin{aligned}
s_k &:= \frac{1}{i^k} \frac{d^k}{dt^k} \ln \mathbb{E}[e^{it\xi_o}] \bigg|_{t=0}, \\
s_k + a_k &:= \frac{1}{i^k} \frac{d^k}{dt^k} \ln \mathbb{E}[e^{it\xi_d}] \bigg|_{t=0}.
\end{aligned}
\end{equation}

\begin{thm}[Landon-Sosoe] \lab{thm:LS}
Fix $\kappa > 0$. Let $H$ be a real symmetric Wigner matrix as above, and $\kappa N \le i \le (1-\kappa)N$. Then,

\[
N \mathbb{E}_H [\lambda_i - \gamma_i] = \frac{1}{2 \pi \rho_{\text{sc}}(\gamma_i)} \arcsin \left( \frac{\gamma_i}{2} \right) - \frac{1}{2 \rho_{\text{sc}}(\gamma_i)} + \frac{s_4}{4} (\gamma_i^3 - 2 \gamma_i) + \frac{a_2 - 1}{2} \gamma_i + o(1)
\]

and

\[
\text{Var}_H (\lambda_i) = \text{Var}_{\text{GOE}} (\lambda_i) + \frac{s_4}{8N^2} \gamma_i^2 + \frac{a_2 - 1}{N^2} + o(N^{-2}).
\]
\end{thm}
Also, $$\E_{GUE} [\la_i] = \gamma_i + n^{-1}C_{i,n} +o(n^{-1}),$$ where $C_{i, n} =  - \frac{1}{2 \rho_{\text{sc}}(\gamma_i)} + \frac{s_4}{4} (\gamma_i^3 - 2 \gamma_i) + \frac{a_2 - 1}{2} \gamma_i$ where $s_4$ equals the sum of the fourth cumulants of the real and imaginary parts of the off-diagonal entries of the matrix.

\end{document}